%% file: main.tex
\UseRawInputEncoding

\documentclass{article}
\pdfoutput=1 

\usepackage[numbers]{natbib}

\usepackage{lineno}
\usepackage{float} 
\usepackage{amssymb}
\modulolinenumbers[5]

\usepackage{amsmath}
\usepackage[ruled, vlined, lined, linesnumbered, commentsnumbered]{algorithm2e} 
\usepackage{multirow}
\usepackage{adjustbox}

\usepackage[figuresleft]{rotating}
\usepackage{caption}
\usepackage{dirtytalk}
\usepackage{authblk}

\makeatletter
\newcommand{\removelatexerror}{\let\@latex@error\@gobble}
\makeatother

\usepackage[margin=1.5in]{geometry}

\newbox{\myorcidaffilbox}
\sbox{\myorcidaffilbox}{\large\includegraphics[height=1.7ex]{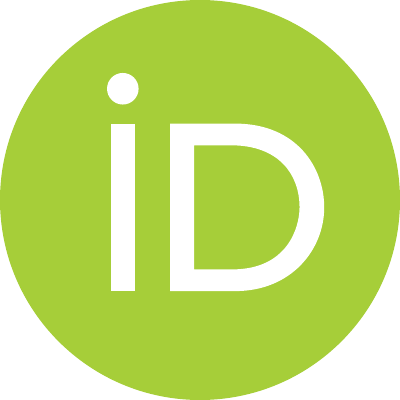}}
\newcommand{\orcidaffil}[1]{%
  \href{https://orcid.org/#1}{\usebox{\myorcidaffilbox}}}
  
\usepackage[utf8]{inputenc}

\usepackage{lineno}
\usepackage{float} 
\usepackage{amssymb}
\modulolinenumbers[5]
\usepackage[ruled, vlined, lined, linesnumbered, commentsnumbered]{algorithm2e} 
\usepackage{multirow}
\usepackage{adjustbox}
\usepackage{color}
\usepackage{amsmath}
\usepackage{dirtytalk}
\usepackage{mathtools}
\usepackage{textcomp}

\usepackage{hyperref}
\hypersetup{pdfauthor={Mahmoud Shaqfa},pdftitle={Disk harmonics for rough surfaces}}
\usepackage[x11names]{xcolor}                     
\usepackage{tikz}
\usetikzlibrary{arrows.meta,calc,patterns,angles,quotes} 
\usepackage{standalone}
\usepackage{amsmath}

\usepackage[normalem]{ulem} 

\newcommand{\rom}[1]{%
  \textup{\uppercase\expandafter{\romannumeral#1}}%
}

\newcommand{\norm}[1]{\left\lVert#1\right\rVert}

\usepackage[figuresleft]{rotating}
\usepackage{caption}

\usepackage{todonotes}

\renewcommand{\vec}[1]{\mathbf{#1}}

\title{Disk Harmonics for Analysing Curved and Flat Self-affine Rough Surfaces and the Topological Reconstruction of Open Surfaces\\~\\\small{[Preprint]}}

\begin{document}

\author{Mahmoud Shaqfa \thanks{ Massachusetts Institute of Technology (MIT), Department of Mechanical Engineering, Cambridge, MA 02139, USA } \orcidaffil{0000-0002-0136-2391} 
Gary P.T. Choi \thanks{The Chinese University of Hong Kong, Department of Mathematics, Hong Kong} \orcidaffil{0000-0001-5407-9111}
Guillaume Anciaux \thanks{\'Ecole Polytechnique F\'ed\'erale de Lausanne (EPFL), Computational Solid Mechanics Laboratory (LSMS), Lausanne, Switzerland}\orcidaffil{0000-0002-9624-5621}
Katrin Beyer \thanks{\'Ecole Polytechnique F\'ed\'erale de Lausanne (EPFL), Earthquake Engineering and Structural Dynamics (EESD), Lausanne, Switzerland}\orcidaffil{0000-0002-6883-5157}
}








\date{}
\maketitle

\input{abstract}



\section{Introduction}
\label{SEC:INTRO}
\input{introduction}

\section{Planar disk parameterisation}
\label{SEC:PARAMETERISATION}
\input{planar_disk_parameterisation}

\section{Fourier-Bessel basis functions}
\label{SEC:FOURIER_BESSEL}
\input{fourier_bessel_basis}

\section{Shape descriptors and their interpretation}
\label{SEC:SHAPE_DESC}
\input{shape_descriptors}



\section{Fractal surfaces: Discussions and validations}
\label{SEC:Fract_RES}
\input{fractal_wavelength}

\section{Surface morphology and reconstruction benchmarks}
\label{SEC:Morph_RES}
\input{morph_results}

\section{Conclusions}
\label{SEC:Conc}
\input{conclusions}


\section{Reproducibility}
\label{SEC:Repro}
\input{reproducibility}

\section*{Acknowledgement}
We would like to thank prof. Ketson R. M. Dos Santos (University of Minnesota) for his thorough revision and feedback on the paper.

\appendix
\newpage
\section{Supplemental derivations}
\label{SEC:APPX}
\input{appendix}


\bibliographystyle{unsrtnat}
\bibliography{main}

\end{document}

%% file: abstract.tex
\begin{abstract}
  When two bodies get into contact, only a small portion of the apparent area is actually involved in producing contact and friction forces, because of the surface roughnesses. It is therefore crucial to accurately describe the morphology of rough surfaces for instance by extracting the fractal dimension and the so-called \textit{Hurst} exponent which is a typical signature of rough surfaces. This can be done using harmonic decomposition, which is easy for periodic and nominally flat surfaces since \textit{Fourier transforms} allow fast and reliable decomposition. Yet, it remains a challenging task in the general curved and non-periodic cases, where more appropriate basis functions must be used. In this work, disk harmonics based on Fourier-Bessel basis functions are employed for decomposing open single-edge genus-0 surfaces (no holes) as a practical and fast alternative to characterise self-affine rough surfaces with the power Fourier-Bessel spectral density. An analytical relationship between the power spectrum density decay and the Hurst exponent is derived through an extension of the Wiener-Khinchin theorem, in the special case where surfaces are assumed self-affine and isotropic. Finally, this approach is demonstrated to successfully measure the fractal dimension, and the \textit{Hurst} exponent, without introducing typical biases coming from basis functions boundary conditions, surface discretisation or curvature of the surface patches. This work opens the path for contact mechanics studies based on the Fourier-Bessel spectral representation of curved and rough surface morphologies. All implementation details for this method are available under GNU LGPLv3 terms and conditions.
\end{abstract}



{\bf Keywords:}
Surface morphology, disk harmonics, self-affine surfaces, Hurst exponent, fractal surfaces, contact mechanics

%% file: introduction.tex
\par
The morphology of natural surfaces is commonly studied using spectral methods, which have recently been exploited in engineering and mechanical applications for studying the behaviour of granular matter, such as sand, gravel and powder-like materials (e.g. \cite{CAPOZZA2021}). However, less attention has been devoted to analysing nominally flat and open surfaces, even though such spectral methods can assist in understanding the frictional behaviour exerted by non-closed surfaces.
\par
In mechanical engineering applications where two rough surfaces are pressed against each other, the resulting friction is known to be highly dependent on surface morphology. The frictional resistance is in reality proportional to the true contact area, which is much smaller than the apparent area \cite{dieterich_imaging_1996, weber_molecular_2018-1, Persson2005} because of the presence of rough asperities. Several models exist to predict the contact-load relation either analytically \cite{persson_theory_2001, barber_thermal_1969, barber_multiscale_2013, greenwood_contact_1966, nayak_random_1971, xu_statistical_2014} or numerically \cite{stanley_fft-based_1997, polonsky_numerical_1999, pei_finite_2005-1, jacq_development_2002, chaise_contact_2011, hyun_elastic_2007, frerot_fourier-accelerated_2019, paggi_coefficient_2010, yastrebov_contact_2012-1, yastrebov_contact_2014, yastrebov_infinitesimal_2015, yastrebov_accurate_2017}. However, these models mostly relied on the hypothesis of a Gaussian, nominally flat surface, with a self-affine geometric characterisation, which manifests with height profiles having a power-law spectral content. Therefore, the spectral decomposition of the surface morphology is crucial when analysing the contact mechanics aspects of rough surfaces which eventually lead to complex deformations such as cracks within concrete or rock and structures built from natural stones.
\par
However, experimentally scanned real surfaces are usually embedding a curvature and are aperiodic. As such, the traditional power spectral density (PSD) method, which is the standard method for analysing surface morphologies, suffers from several limiting drawbacks \cite{Jacobs2017}:
\begin{itemize}
    \item The employed Fourier transformation assumes that the signal is periodic; when applying the method to aperiodic real surfaces, the required periodicity must be imposed on the data using a windowing function, see Prabhu (2018) \cite{Prabhu2018} and Elson \& Bennett (1995) \cite{Elson1995}. These windowing functions \say{taper} the data close to the boundaries of the image and work like a low-pass filter. Studies indicate that the choice of windowing function affects the PSD results \cite{Jacobs2017}.
    \item For the PSD method, the tilt of the input signal must be removed in a pre-processing step. Jacobs et al. (2017) \cite{Jacobs2017} proposed a method for removing this tilt by fitting and subtracting a plane from the signal to adjust the average tilt of the processed surface to zero. However, such corrections are difficult to apply to curved surfaces and render the PSD results sensitive to the pre-processing method.
\end{itemize}
\par
Next to the PSD method for studying 3D rough surfaces, spherical and hemi-spherical harmonics analyses have also been proposed. However, these basis functions also come with their own setbacks (reviewed in Shaqfa et al. [2021] \cite{shaqfa2021b}):
\begin{itemize}
    \item Because spherical harmonics are used for studying closed genus-0 objects, open surfaces cannot be expanded regionally on a closed sphere without proper windowing functions. These windowing functions will alter the actual input signals (surface data) similarly to the traditional PSD when forcing periodic boundary conditions.
    
    \item The convergence of these functions is slow and could be sped up considerably with a method that analyses the morphology at any starting scale from macro to micro and across several scales.
    
    \item Although recent papers employed hemispherical harmonics (variant basis functions from spherical harmonics), these functions cannot accurately represent the morphology of rough surface patches that are nominally flat because a significant distortion is introduced when mapping a nominally flat surface to a hemisphere.
\end{itemize}
\par
To overcome these issues, we recently proposed the use of spherical cap harmonic analysis (SCHA) for treating rough surfaces that are open and nominally flat in Shaqfa et al. (2021) \cite{shaqfa2021b}. In that paper, we showed how the surface analysis and reconstruction accuracy depends on the chosen parametrisation algorithm. That paper also connected the analysis degrees with a physical wavelength of the surface and provided an empirical relationship between the shape descriptors from SCHA and the fractal dimension of the surface.
\par
In Shaqfa et al. (2021) \cite{shaqfa2021b}, we used shallow spherical caps to minimise the distortion between the spacial and parameterisation domains. We treated the determination of the optimal half-angle $\theta_c$ of the spherical cap as an optimisation problem. In that paper, introducing $\theta_c$ as a free parameter provided a flexible mapping approach, though the fact that $\theta_c$ was always determined to be close to zero for flat surfaces means that this method worked only near a singular solution for the spherical basis function. Practically, this indicates that the nominally flat surface should instead be parameterised onto a planar disk instead of a curved shallow spherical cap. Moreover, the additional projection from a topological disk to a spherical cap always added distortion to both the areas and angles of the mesh elements, and the use of Gaussian hypergeometric function $_2F_1(a, b; c; z)$ for the fractional associated Legendre functions made the basis convergence slow and numerically unstable for high degrees.
\par
\begin{figure}[t!]
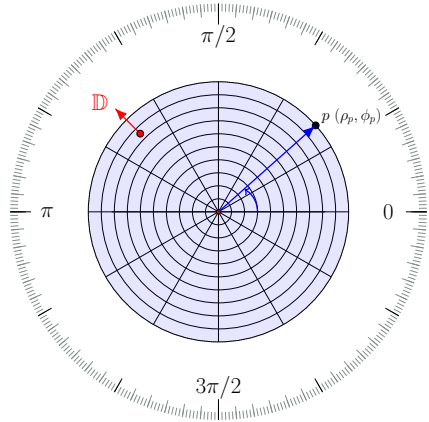

\centering
  \includestandalone[width=0.4\textwidth]{planar_disk_tikz}%
  \caption{A planar unit disk $\mathbb{D}$ \big($\rho(\phi) \leq 1 ~\forall \phi \in [0, 2\pi]$\big) that is centred about the pole in a polar coordinate system. The angle $\phi$ is measured counter-clockwise from the polar axis in radian.}
  \label{FIG:planar_disk}
\end{figure}

To overcome these setbacks, we herein propose the use of disk harmonics (DH), which consist of Fourier-Bessel basis functions over a planar disk in the polar coordinate system (Fig.~\ref{FIG:planar_disk}). These Fourier-Bessel functions are the natural extension of Fourier analysis in the polar coordinate system and are unconditionally stable for any input which solves the stability problem encountered for SCHA \cite{shaqfa2021b}. Also, the mathematical connection between the spherical harmonics and disk harmonics bases comes from assuming small $\theta_c$, the near-pole areas, as they tend to become flat surfaces, then, the associated Legendre functions can be approximated by Bessel functions for large expansion degrees (see Eq.~[9.1.71] in the handbook of Abramowitz and Stegun (1964) \cite{abramowitz_stegun_1964}). Moreover, this method simplifies the parametrisation process by avoiding the extra optimisation step used in SCHA, which maps the planar disk onto the spherical cap using prescribed half-angle $\theta_c$. This paper depends heavily on the concepts introduced in our previous work \cite{shaqfa2021b}, and for this reason, these concepts will be briefly revisited in the current work. To the authors' knowledge, this is the first time these basis functions are used for studying the morphology of open surfaces.
\par
Literature treating an open-surface morphology is scarce. Two main studies by Huang et al. (2006) \cite{REF:6} and Giri et al. (2020) \cite{REF:7} analysed open surfaces using hemispherical harmonics functions, while we treated open and nominally flat surfaces through SCHA \cite{shaqfa2021b}. For more details, we refer the reader to our previous work \cite{shaqfa2021b}, which compares this treatment with other existing methods.
\par
The basis functions used in this paper are the DH and are also known as the drumhead harmonic basis. Though these basis functions have been used in many engineering and physical applications, to our knowledge, it has not been employed for studying the morphology of surfaces. In structural engineering, DH has been used to study the vibration modal shapes of circular membranes, such as in the work by Kirchhoff (1850) \cite{Kirchhoff1850}, Rayleigh (1877) \cite{Rayleigh1877}, Timoshenko (1929) \cite{Timoshenko1929}, Wah (1962) and \cite{Wah1962}. Similar spectral methods were also used for simulating 3D incompressible flow in cylindrical geometries \cite{Blackburn2004}.
\par
DH has also been used in image analysis. For pattern recognition of 2D images, Verrall et al. (1998) \cite{Verrall:98} studied the reconstruction and invariant descriptors and compared the results with other similar orthogonal bases, including Zernike and pseudo-Zernike. Similar approaches also exist that use different orthogonal variants of a radial basis on a disk (\cite{Sheng1994, Ping2002, Ren2003}). The invariant descriptors extracted from DH can also be used beyond image compression and reconstruction. For example, Sun et al. (2020) \cite{Sun2020} used shape descriptors to rigorously generalise convolutional neural networks (CNNs) for image analysis while maintaining affordable computational efforts.
\par
In the following section, Section \ref{SEC:PARAMETERISATION}, we introduce our proposed combination of density-equalising and quasi-conformal maps for achieving bijective disk area-preserving parameterisations. Section \ref{SEC:FOURIER_BESSEL} explains the Fourier-Bessel basis functions as well as the imposed boundary conditions and the corresponding eigenvalues. Section \ref{SEC:EXTRACTION} explains the extraction of the expansion coefficients. In Section \ref{SEC:SHAPE_DESC}, we interpret the computed shape descriptors and discuss the normalisation approach. In Sections \ref{SEC:Fract_RES} and \ref{SEC:Morph_RES}, we validate the method against fractal self-affine surfaces and benchmark topology problems. Sections \ref{SEC:Conc} summarise our conclusions and set our resolutions for future works, respectively. Two additional appendices support our derivations and include miscellaneous results.

%% file: planar_disk_tikz.tex
\begin{tikzpicture}
\foreach \r in {0.5, 1.0,...,5}
  \draw[black, thick] (0,0) circle (\r);

\foreach \a in {0, 1,...,359}
  \draw[Azure4] (\a:7.7) -- (\a:8);
  
\foreach \a in {0, 5,...,359}
  \draw[Azure4] (\a:7.5) -- (\a:8);  

\foreach \a in {0, 30,...,359}
  \draw[thick, black] (\a:7.4) -- (\a:8);  

\foreach \a in {0, 90,...,359}
  \draw[thick, black] (\a:7.2) -- (\a:8);  

\foreach \a in {0, 30,...,359}
  \draw[thick,black] (0, 0) -- (\a:5);

\draw[fill=red] (0,0) circle(0.7mm);

\draw[dotted, draw=black, fill=blue, fill opacity=0.1] (0,0) circle(5);

\foreach \ang/\lab/\dir in {
  0/0/right,
  2/{\pi/2}/above,
  4/{\pi}/left,
  6/{3\pi/2}/below} {
  \draw (0,0) -- (\ang * 180 / 4:4.1);
  \node at (\ang * 180 / 4:6.2) [\dir] {\Huge $\lab$};
}

\node[red] at (-180 / 4.2:-6.2) {\Huge $\mathbb{D}$};
\draw[-{Latex[length=5mm]}, thick, red, line width=0.5mm] (-45:-4.25) -- (-4, 4);
\draw[fill=red] (0,0) circle(0.8mm);
\draw[fill=red] (-45:-4.25) circle(1.3mm);

\draw[fill=black] (2500 / 60:5) circle(1.3mm);
\node[black] at (2500 / 70:6.3) {\LARGE $p~(\rho_p, \phi_p)$};
\draw[-{Latex[length=5mm]}, thick, blue, line width=0.5mm] (0, 0) -- (2500 / 60:5);
\draw[ultra thick, blue, line width=0.5mm, ->] (1.5,0) arc (0:60:1);

\end{tikzpicture}

%% file: planar_disk_parameterisation.tex
In this section, we combine the density-equalising map (DEM) method~\cite{choi2018density} and the quasi-conformal mapping method~\cite{choi2015flash} to put forward a planar disk area-preserving parameterization $\varphi:\mathcal{S} \to \mathbb{D}$, where $\mathcal{S}$ is a simply connected open surface in the form of a triangular mesh, and where $\mathbb{D} = \left\{ \rho, \phi \mid \rho \in [0, 1], \phi \in [0, 2\pi]\right\}$ is the unit disk. In this work the rotation angle $\phi$ is measured about the origin (pole) and counter-clockwise from the polar axis, as shown in Fig.~\ref{FIG:planar_disk}.
The DEM method aims to smoothly map a given surface onto a planar domain based on density diffusion. Suppose $\varrho > 0$ is a density field prescribed on $\mathcal{S}$.
The DEM method finds a parameterization of $\mathcal{S}$ onto the plane such that its Jacobian determinant is proportional to $\varrho$.
In other words, it flattens $\mathcal{S}$ onto a planar domain and deforms it to achieve a prescribed area change.
In particular, to map $\mathcal{S}$ onto the unit disk, we first compute an initial flattening map $g:\mathcal{S} \to \mathbb{D}$ using Tutte embedding~\cite{tutte1963draw}.
We then apply the DEM method to find a mapping $h:\mathbb{D} \to \mathbb{D}$ by iteratively solving the following diffusion equation with the Neumann boundary condition~\cite{choi2018density}:
\begin{equation}
\left\{
\begin{array}{ll}
     \frac{\partial{\varrho}}{\partial t} = \Delta \varrho &  \text{ on } \ \ g(\mathcal{S}\setminus \partial \mathcal{S}), \\
     \nabla \varrho \cdot \hat{\mathbf{n}} = 0 & \text{ on } \ \ g(\partial \mathcal{S}),
\end{array}\right.
\end{equation}
where $\hat{\mathbf{n}}$ is the unit outward normal. Throughout the iterations, the updated density field can also update the vertex positions $\mathbf{x}$ on the planar domain over time~\cite{gastner2004diffusion}:
\begin{equation}
    \mathbf{x}(t) = \mathbf{x}(0) - \int_0^t \frac{\nabla \varrho}{\varrho} d\tau,
\end{equation}
where $\mathbf{x}(0)$ denotes the initial flattened domain. As $t \to \infty$, $\varrho$ is equalised, resulting in the density-equalising map $h(\mathbf{x}(0)) = \mathbf{x}(\infty)$. In particular, with the initial density set using the area of the triangular faces: $\varrho(T) = \frac{\text{Area}(T)}{\text{Area}(g(T))}$ for each triangle $T$ on $\mathcal{S}$, the planar domain will be deformed based on the local area of $\mathcal{S}$, thereby resulting in a smooth area-preserving map $h \circ g: \mathcal{S} \to \mathbb{D}$~\cite{choi2018density,choi2020area}. 

As mentioned in Choi et al.~\cite{choi2018density}, there is no theoretical guarantee about the bijectivity of the density-equalising map $h$. To ensure the bijectivity of $h$, we follow an idea presented in an earlier work of Choi at al.~\cite{choi2015flash} and compute the Beltrami coefficient $\mu_h:\mathbb{D} \to \mathbb{C}$ of the mapping $h$, which is a complex-valued function that measures the quasi-conformal distortion caused by mapping $h$. In particular, $\|\mu_h\|_{\infty} < 1$ if and only if $h$ is a bijective map~\cite{lehto1973quasiconformal}. Therefore, we herein propose a simple method for achieving a bijective density-equalising map by first rescaling all large $|\mu_h(z)|$ to a smaller magnitude to update the Beltrami coefficient $\tilde{\mu}_h$ and then applying the linear Beltrami solver (LBS)~\cite{choi2015flash} to update the density-equalising map $\tilde{h}:\mathbb{D} \to \mathbb{D}$ associated with $\tilde{\mu}_h$. Note that this procedure only locally corrects the non-bijectivity and does not alter the overall density-equalising effect of $h$. 

\begin{figure}
\centering
\includegraphics[width=0.9\textwidth]{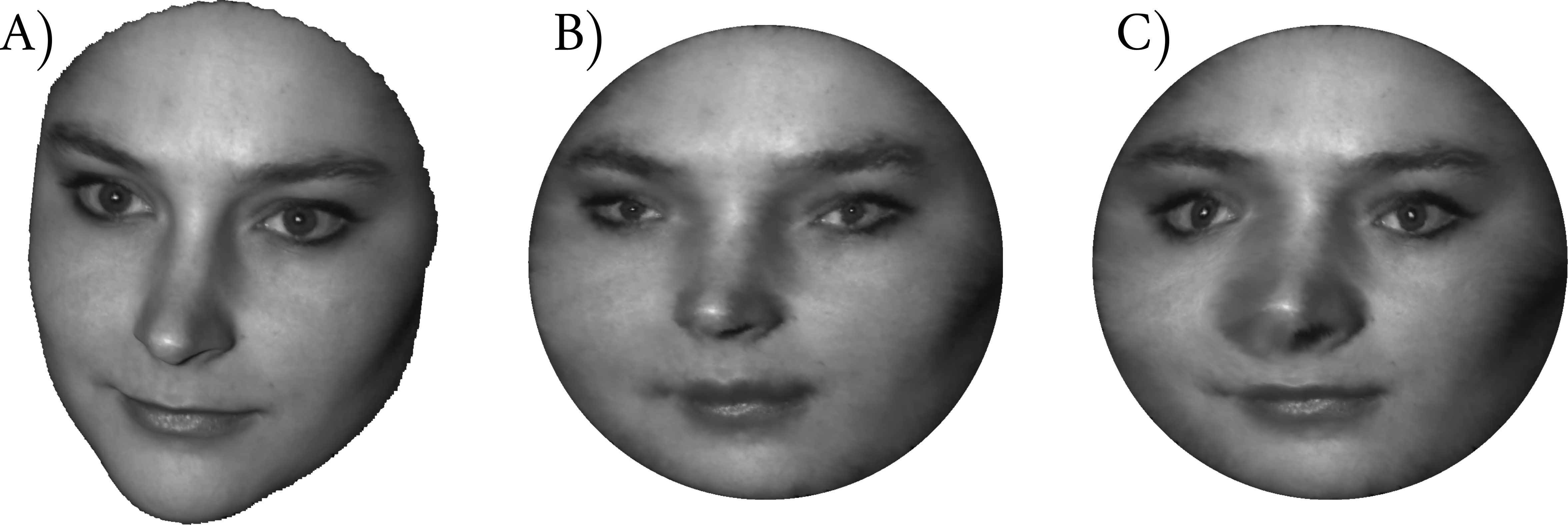}%
\caption{A planar disk area-preserving parameterisation $\varphi:\mathcal{S} \to \mathbb{D}$ of a benchmark face named Sophie (the mesh was retrieved from \cite{Riemannmapper}). A) The input surface mesh $\mathcal{S}$. B) The initial flattening map $g:\mathcal{S} \to \mathbb{D}$. C) The final area-preserving map $\varphi:\mathcal{S} \to \mathbb{D}$ given by $\varphi = \tilde{h} \circ g$, where $\tilde{h}$ is the updated map after enforcing the bijectivity of the density-equalising map $h$.
\label{FIG:sophie_benchmark}
}
\end{figure}

Finally, the desired parameterisation $\varphi:\mathcal{S} \to \mathbb{D}$ is given by $\varphi = \tilde{h} \circ g$, where the choice of the density in the DEM method produces $\varphi$ as a smooth area-preserving mapping of $\mathcal{S}$ onto $\mathbb{D}$. Also, by Tutte's spring theorem~\cite{tutte1963draw} and quasi-conformal theory~\cite{lehto1973quasiconformal}, both $g$ and $\tilde{h}$ are guaranteed to be bijective, such that $\varphi$ is also bijective. Figure \ref{FIG:sophie_benchmark} shows the parameterisation steps of a benchmark triangular mesh.

%% file: fourier_bessel_basis.tex
\par
In this section, the mathematical derivation of the DH basis functions will be briefly described. For a more extensive derivation of these previously established basis functions, we refer the reader to the following textbooks: Tolstov and Silverman (1976) \cite{REF:1}, Spiegel (1974) \cite{REF:2}, Brown and Churchill (1993) \cite{REF:3} and Wolf (1979) \cite{REF:4}. Many variants for basis functions of flat disks can be found in the literature, such as the Zernike and Chebyshev basis; comparisons can be found in Boyd and Yu~\cite{Boyd2011}. In this paper, the DH analysis will be made with the Fourier-Bessel basis functions $D_m^k (\rho, \phi): \mathbb{D} \rightarrow \mathbb{R}$ which are bijections expressed with polar coordinates on the unit disk, therefore taking the form:
\begin{equation}
  D_m^k (\rho, \phi) = \underbrace{N_m^k}_{\text{Normalisation}} \overbrace{J_m\big( l(m)_k~\rho \big)}^{\text{Bessel}} ~\underbrace{e^{i m \phi}}_{\text{Fourier}}.
  \label{EQN:FOUR_BESS_final}
\end{equation}
Where $l(m)_k$ is a real number which is associated with order $m$ and degree $k$ of the basis function, as will be explained below. 
These basis functions allow a decomposition of any 2D function (or signal) parametrised in polar coordinates using the following series expansion:
\begin{equation}
    f(\rho, \phi) = \sum_{k = 0}^{\infty} \sum_{m = -k}^{k} q_m^k~ D^{k}_{m} (\rho, \phi),
\end{equation}
where $q_m^k$ represents the expansion coefficients for different degrees and orders, as we will explain next.
\par
The basis functions are classically chosen as eigenfunctions of the Laplacian operator, which can be found by solving $\nabla^2 f(\rho, \phi) = - m^2f$ (Helmholtz equation) and is expressed in polar coordinates as follows:
\begin{equation}
  \frac{\partial^2{f}}{\partial{\rho^2}} + \frac{1}{\rho}\frac{\partial{f}}{\partial{\rho}} + \frac{1}{\rho^2}\frac{\partial^2{f}}{\partial{\phi^2}} = -m^2 f,
    \label{EQN:LAPLAC}
\end{equation}
where $m$ is the basis order. To ensure periodicity and continuity in the polar coordinates, viable eigenfunctions must also satisfy the following periodic boundary conditions:
\begin{equation}
    f(\rho, \phi) = f(\rho, \phi + 2\pi) \qquad \text{and} \qquad     \frac{\partial f(\rho, \phi)}{\partial \phi} = \frac{\partial f(\rho, \phi + 2\pi)}{\partial \phi}.
    \label{EQN:BC1}
\end{equation}
The eigenfunctions satisfying Eq.~(\ref{EQN:LAPLAC}) with the associated boundary conditions are called harmonic functions and can be found by separating variables $\rho$ and $\phi$ such that $f(\rho, \phi) = R(\rho)\Phi(\phi)$. The separation gives the following two differential equations:
\begin{align}
  &\Phi^{\prime\prime}(\phi) - m^2 \Phi(\phi) = 0 \label{EQN:euler}\\
  &\rho^2 R^{\prime\prime}(\rho) + \rho R^{\prime}(\rho) + (\rho^2 - m^2) R(\rho) = 0, \label{EQN:SL}
\end{align}
where the second one is the \textit{Sturm-Liouville} problem.
\par
Equation (\ref{EQN:euler}) results in Euler's formula, $\Phi(\phi) = e^{i m \phi}$, where $m$ becomes a positive integer when periodic boundary conditions are enforced. Thus, $m$ can be interpreted similarly to the case of spherical harmonics\cite{shaqfa2021b}. Bessel functions of the first and second kind \cite{REF:1} are solutions of the \textit{Sturm-Liouville} problem in polar coordinates. However, the basis functions must be finite (non-singular) at the origin ($\rho = 0$), and this can only be satisfied by the Bessel function of the first kind (see Boyd and Yu~\cite{Boyd2011}). The general solution of the \textit{Sturm-Liouville} problem is, therefore:
\begin{equation}
    R(\rho) = J_m \big(l(m)~\rho \big),
    \label{EQN:LAP_BESSEL}
\end{equation}
where $J_m$ is the Bessel function of the first kind, and $l(m)$ is a constant associated with order $m$ which will be determined while imposing other boundary conditions. Indeed, obtaining a solution on the unit disk requires boundary conditions to be imposed at the edge of the unit disk, i.e. where $\rho = 1$. As in \cite{shaqfa2021b}, the Neumann boundary condition enforcing vanishing slopes on the disk edge is now imposed. This allows us to determine all possible values for the eigenvalue $l(m)$ by finding the roots of the following equation:
\begin{equation}
    J_m^{\prime} \big(l(m) \cdot 1 \big) = 0.
    \label{EQN:BC4}
\end{equation}
This leads to the formal definition of $l(m)_k$ as the $k$-the root for the derivative of the $m$-th solution to the \textit{Sturm-Liouville} problem.
\par
The functions $J_m(x)$ and $J_m^{\prime}(x)$ have an infinite number of roots, and for large values of $x$, the roots are approximately separated by $\pi$ (see Eq.~[\ref{EQN:Bessel_ASYM}] in \ref{APP:wavelengths} for more details). Tabulated data for the roots of $J_m^{\prime}(x)$ can be found in the literature\cite{REF:1, Morgenthaler1963}. These roots $l(m)_k$ can be computed with a simple stepping solver based on Mueller's method \cite{REF:66}, which is similar to the one used in Shaqfa et al.~\cite{shaqfa2021b}.
\par
The harmonic functions described so far still need to be normalised by the factors $N_m^k$ to obtain the desired orthonormal basis of functions. With the classical integral scalar product defined as
\begin{equation}
\left< f, g\right> = \int_0^1  \int_0^{2\pi} f(\rho, \phi) \overline{g}(\rho, \phi) \rho d\theta d\rho,
\end{equation}
where, $\overline{g}$ is the complex conjugate of $g$. The norm of the eigenfunctions can therefore be split to become:
\begin{equation}
|| J_m(l(m)_k \rho) e^{im\phi}||^2 = 2\pi \int_0^1  \left[J_m(l(m)_k \rho)\right]^2  \rho d\rho. 
\end{equation}
Bessel functions' identities, as employed in \cite{REF:1, REF:4}, allow obtaining:
\begin{align}
  \int_0^1  \left[J_m(l(m)_k \rho)\right]^2  \rho d\rho &= \frac{1}{2}\left\{\underbrace{\left[J_{m}^{\prime}\big(l(m)_k\big)\right]^2}_{= 0} + \Bigg(1 - \frac{m^2}{l(m)_k^2}\Bigg) \left[J_m \big(l(m)_k\big)\right]^2\right\}.\\
 &= \frac{1}{2}\Bigg(1 - \frac{m^2}{l(m)_k^2}\Bigg) \left[J_m \big(l(m)_k\big)\right]^2.
    \label{EQN:NORM_BESSEL}
\end{align}
where the application of the Neumann boundary condition and Eq.~(\ref{EQN:BC4}) brought the simplification. The final normalisation factor becomes:
\begin{equation}
    N_m^k = \frac{J_m^{-1} \big(l(m)_k\big)}{\sqrt{\pi \Big(1 - \frac{m^2}{l(m)_k^2}\Big)}},
    \label{EQN:NORM}
\end{equation}
which gives the orthonormal Fourier-Bessel harmonic functions as defined in
Eq.~(\ref{EQN:FOUR_BESS_final}). 
\par
\begin{figure*}[!ht]
    \centering
    \includegraphics[width=0.8\textwidth]{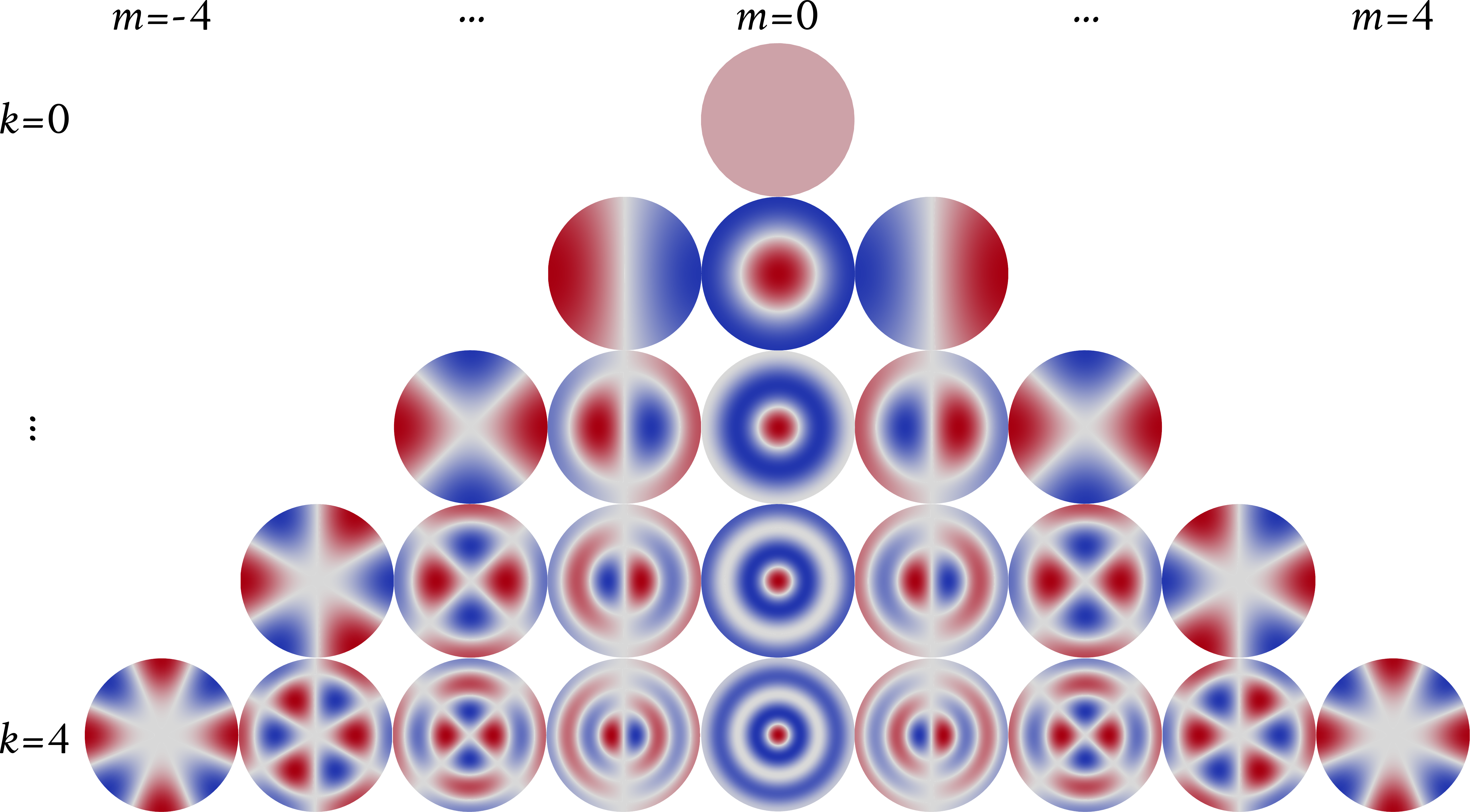}
    \caption{The normalised Fourier-Bessel basis functions $\mathfrak{R}\Big(D^{k}_{m} (\rho, \phi) \Big)$.}
    \label{FIG:BESSEL_FOURIER_BASIS}
\end{figure*}
\begin{figure*}[!ht]
    \centering
    \includegraphics[width=0.9\textwidth]{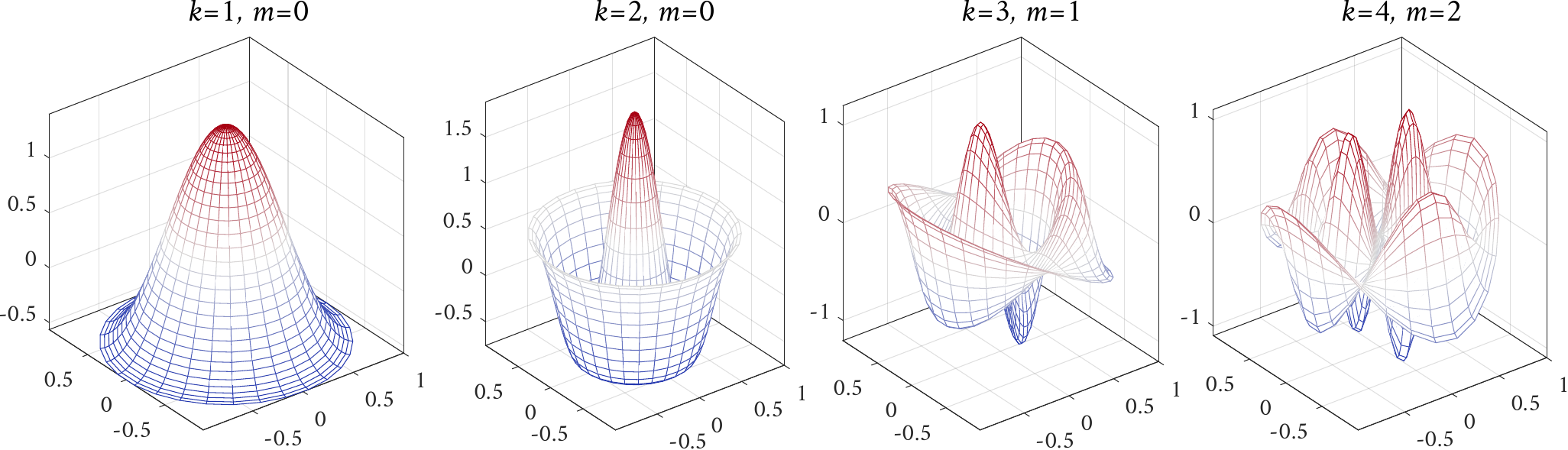}
    \caption{The normalised Fourier-Bessel basis functions $\mathfrak{R}\Big(D^{k}_{m} (\rho, \phi) \Big)$. The complex part of the base functions is identical to the real ones except that they are rotated by $\pi/(2|m|)$ about the $z$-axis.}
    \label{FIG:BESSEL_FOURIER_BASIS_3D}
\end{figure*}

Finally, the Condon-Shortley identity allows relating the
basis functions of order $m$ and $-m$:
\begin{equation}
    D_{-m}^k (\rho, \phi) = (-1)^m~ \overline{D}_{|m|}^{k} (\rho, \phi),
    \label{EQN:FOUR_BESS_NEG_ORDR}
\end{equation}
where $\overline{D}_{|m|}^{k}$ is the complex conjugate of $D_{|m|}^{k}$ the basis function, which
implies that only positive order functions have to be numerically extracted and stored. Some Fourier-Bessel basis functions are shown in Fig. \ref{FIG:BESSEL_FOURIER_BASIS}, and the 3D profile of a selected set of four such functions is plotted in Fig. \ref{FIG:BESSEL_FOURIER_BASIS_3D}.

\section{Extraction of the expansion coefficients} \label{SEC:EXTRACTION}
\par
A bi-directional signal $\vec{f}: \mathbb{D} \rightarrow \mathbb{R}^3$ that is parameterised onto a unit disk in the polar coordinates can then be represented by the series expansion:
\begin{equation}
    \vec{f}(\rho, \phi) = \sum_{k = 0}^{\infty} \sum_{m = -k}^{k} \vec{q}^k_m~ D^{k}_{m} (\rho, \phi).
    \label{EQN:SIGNAL_INF_SUM}
  \end{equation}
where every coefficient $\vec{q}^k_m$ is a vector in $\mathbb{R}^3$ (hence the bold notation).
The signal $\vec{f}(\rho, \phi)$ is parameterised onto 2D polar coordinates over a planar unit disk, such that:
\begin{equation}
    \vec{f}(\rho, \phi) = 
        \begin{pmatrix}
        f_x(\rho, \phi)\\
        f_y(\rho, \phi)\\
        f_z(\rho, \phi)
        \end{pmatrix}.
    \label{EQN:CART_VECTOR}
  \end{equation}
In order to determine the expansion coefficients, each direction will be treated separately in what follows. The approximation of Eq.~(\ref{EQN:SIGNAL_INF_SUM}) as a truncated series with a maximum index $k = K_{max}$ leads to:
\begin{equation}
    f(\rho, \phi) = \sum_{k = 0}^{K_{max}} \sum_{m = -k}^{k} q_{m}^k~ D_{m}^{k} (\rho, \phi) = D(\rho, \phi) Q.
    \label{EQN:SIGNAL}
\end{equation}
Here, $q_{m}^k$ are the expansion coefficients, which can be organised in a vector $Q$ so that $q_m^k = Q_{k^2 + k + m}$ as described in 
\cite{REF:5}. Because the basis functions are orthonormal, these coefficients can be computed from the inner dot product:
\begin{equation}
    q_{m}^k =~ \langle f(\rho, \phi),  D^{k}_{m} (\rho, \phi)\rangle,
    \label{EQN:COEF_INTG1}
    \qquad \text{with} \qquad
     \langle f, g \rangle = \int_0^{2\pi} \int_0^1 f(\rho, \phi) ~ g(\rho, \phi)\rho~ d \rho~ d \phi.
\end{equation}
Similar expansions can be obtained for the other directions $y$ and $z$.

However, the discretization of the planar disk resulting from the previously described parameterisation algorithm does not form a uniform distribution of the vertices, such that orthogonality is not guaranteed. Instead, the optimal coefficients can be computed in the least-square sense, which provides a statistical regularisation. The squared distance is defined using the 2-norm $\norm{f}^2 = \langle f, f \rangle$ \cite{REF:3}:
\begin{equation}
\Gamma(Q) = \sum_{i \in \mathcal{V}} \norm{f(\rho_i, \phi_i) - D(\rho_i, \phi_i) Q}^2,
\end{equation}
which leads to the optimal coefficients $Q$ when minimising $\Gamma$. This further leads to a linear system to solve:
\begin{equation}
\forall j \in \mathcal{V} \qquad 2\sum_{i \in \mathcal{V}} D^T(\rho_j, \phi_j) \cdot \Big(f(\rho_i, \phi_i) - D(\rho_i, \phi_i)Q\Big) = 0.
\end{equation}
Taking $f(\rho, \phi)$ as $f_x(\rho, \phi)$, $f_y(\rho, \phi)$ and $f_z(\rho, \phi)$ brings $(K_{\max}+1)^2$ unknown coefficients per direction in the space, or in other words, brings one set of coefficients for each $x$, $y$ and $z$. This 
leads to the solution expressed with matrices as described in Brechb{\"u}hler et al. \cite{REF:5}:
\begin{equation}
    B^{T} B Q = B^{T} V,
    \label{EQN:LSF}
\end{equation}
where $Q$ holds the desired coefficients, while $V$ is the matrix of the coordinates of the $n_v$ vertices and $B$ is the matrix holding the DH basis functions:
\begin{equation}
    V = \begin{pmatrix}
        x_0    & y_0    & z_0\\ 
        x_1    & y_1    & z_1\\ 
        \vdots & \vdots & \vdots\\
        x_{n_{v}-1} & y_{n_{v}-1} & z_{n_{v}-1}
    \end{pmatrix},
    \label{EQN:VERTS}
\quad 
    B = \begin{pmatrix}
    D(\rho_0, \phi_0)\\
    D(\rho_1, \phi_1)\\
    \vdots \\
    D(\rho_{n_v-1}, \phi_{n_v-1})\\    
    \end{pmatrix},
\quad
    Q = 
    \begin{pmatrix}
        q^0_{x}    & q^0_{y}    & q^0_{z}\\ 
        q^1_{x}    & q^1_{y}    & q^1_{z}\\ 
        \vdots & \vdots & \vdots\\
        q^{j}_{x} & q^{j}_{y} & q^{j}_{z}
    \end{pmatrix}.
\end{equation}
 where $x_i = f_x(\rho_i, \phi_i)$, $y_i = f_y(\rho_i, \phi_i)$, and $z_i = f_z(\rho_i, \phi_i)$.

%% file: shape_descriptors.tex
\par
The DH analysis produces $3(K_{\max}+1)^2$ coefficients $q^k_m$, which are referred to as Fourier weights or amplitudes or more specifically polar Fourier descriptors (PFD). These weights can be used to construct invariant shape descriptors that characterise the surface similarly to coefficients obtained by spherical \cite{REF:5} and spherical cap \cite{shaqfa2021b} harmonics analyses.
\par
\begin{figure}
\centering
\includegraphics[width=0.9\textwidth]{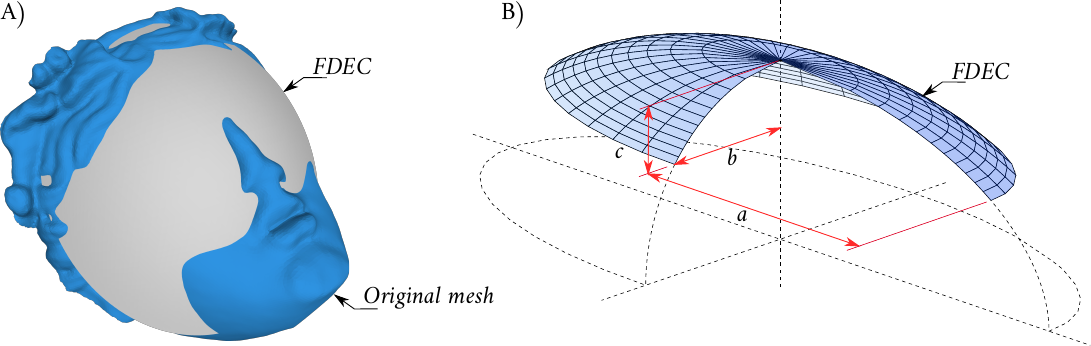}%
\caption{Fitting the first degree ellipsoidal cap (FDEC) defined by an expansion at $k = 1$. A) We show how the FDEC, reconstructed at $k=1$, represents an ellipsoidal cap fitted over the input surface. B) The size of the cap is determined by $a$, representing the half-major axis, $b$, representing the half-minor axis, and $c$, representing the ellipsoidal cap depth, where $|c| \leq |b| \leq |a|$. Inset (B) was retrieved from Shaqfa et al.~\cite{shaqfa2021b}.
\label{FIG:FDEC}
}
\end{figure}

First, the zero-degree weights ($q_{0,x}^0$, $q_{0,y}^0$ and $q_{0,z}^0$) define a constant shift of the geometric centre of the surface with regard to the origin of the coordinate system. Then, the weights for $k = 1$ for the three directions carry curvatures characteristic since they allow to define an ellipsoidal cap, here termed the first degree ellipsoidal cap (FDEC) which follows the equation:
\begin{equation}
  \vec{f}(\rho, \phi) = \vec{q}_{-1}^1 D_{-1}^1(\rho, \phi) + \vec{q}_0^1 D_0^1(\rho, \phi) + \vec{q}_1^1 D_1^1(\rho, \phi). \label{eq:fdec-eq}
\end{equation}
Figure \ref{FIG:FDEC} illustrates such a cap, which has already been introduced in Shaqfa et al. (2021) \cite{shaqfa2021b}. It is visible on this figure that the global curvature is captured whereas the details of the original mesh have been discarded. The so-called FDEC sizes are the lengths of the ellipsoidal cap half-axes. These can be extracted by either manipulation of Eq.~(\ref{eq:fdec-eq}) (for details, see \ref{APP:FDEC_size} and \cite{shaqfa2021b, REF:5, REF:8}) or direct measurement of reconstruction at $k=1$. Following a similar approach to \cite{shaqfa2021b}, based on the identified weights at $k = 1$, leads to extracting the eigenvalues of the matrix $\vec{A}^t\vec{A}$ with:
\begin{equation}
    \label{EQN:FDE_MAT}
    \vec{A} = \frac{1}{2}N^1_1 l_1(1) \left(
\begin{matrix}
  \vec{q}_{-1}^{1}-\vec{q}^{1}_{1} \\
  i(\vec{q}_{-1}^{1} + \vec{q}^{1}_{1})
\end{matrix}
  \right).
\end{equation}
The matrix $\vec{A}$ is defined above with two line vectors and allows to extract
the directions of the two principal axes tangent to the ellipsoidal cap, with
their lengths being the square root of the eigenvalues $a^2 = |\lambda_1| \geq b^2 = |\lambda_2|$. The last eigenvalue (the depth of the cap) is found to be $c^2 = 2N_0^1l(0)_1^2$. This approach, as demonstrated in \ref{APP:FDEC_size}, does not provide good estimates via the disk harmonics bases. Alternatively, as discussed in \ref{APP:FDEC_size},
we can use the oriented bounding box (OBB) method via, for example, the principal component analysis (PCA) for fitting n-dimensional ellipsoids that can result in the size of the cap after reconstructing the mesh at $k=1$.
\par
Scalar shape descriptors can be derived from the expansion coefficients in $Q$, which are invariant to rotations. These are usually taken along each axis (for details, see Shaqfa et al.~\cite{shaqfa2021b}):
\begin{equation}
\forall i \in \{x, y, z\}\qquad \hat{D}_{k, i} = \sqrt{\sum_{m = -k}^{k} ||q^k_{m,i}||^{2}}.
\end{equation}
but scalar amplitudes can also be simply extracted:
\begin{equation}
    \hat{D}_{k}^2 = \hat{D}_{k, x}^2 + \hat{D}_{k, y}^2 + \hat{D}_{k, z}^2.
    \label{EQN:DESCRIPTORS}
\end{equation}
However, these descriptors in Eq.~(\ref{EQN:DESCRIPTORS}) are influenced by the curvature. Normalising them is particularly useful in morphological studies comparing various surfaces at different scales. However, it is not generally used in tribology applications, for instance to compute the power spectral density (PSD) while compensating for the curvature. To account for curvature, the normalisation is simply applied separately in each directions, as explained in \ref{APP:analytic_cap}, leading to the modified shape descriptor:
\begin{equation} \label{EQN:SHAPE_DESC}
    D_{k}^2 = \frac{\hat{D}_{k, x}^2}{\hat{D}_{1, x}^2} + \frac{\hat{D}_{k, y}^2}{{\hat{D}_{1, y}^2}} + \frac{\hat{D}_{k, z}^2}{{\hat{D}_{1, z}^2}}, ~\forall k>1.
\end{equation}
This expression deviates from what we normally find in the literature to account for size normalisation, as we propose to normalise each axis separately along the principal axes of the FDEC. In the special case where the cap radii are the same in all directions, Eq.~(\ref{EQN:SHAPE_DESC}) gives the following relation:
\begin{equation}
    {D}_{k}^2 = \kappa(\hat{D}_{k, x}^2 + \hat{D}_{k, y}^2 + \hat{D}_{k, z}^2),
    \label{EQN:NORM_CAP_DESC}
\end{equation}
where $\kappa = 1/R^2_{\theta_c}$ is the Gaussian curvature of the spherical surface.

%% file: fractal_wavelength.tex
\par
For rough surfaces that are nominally flat, the average statistical properties can be determined using height-height autocorrelation functions (ACF) defined over the spatial domain of the surfaces. The Wiener-Khintchine theorem relates the power spectrum of random surfaces to its corresponding ACF. Persson et al.~\cite{Persson2005} showed that the radial decay of the PSD for an isotropic surface is a function of the Hurst exponent.
\par
In this section, we define a similar relationship between the PSD defined from the Fourier-Bessel basis functions and the Hurst exponent. We then analyse example surfaces and extract their fractal dimensions (Hurst exponent) via the previously described decomposition, with flat surfaces at first. We also explore the impact on the Hurst exponent of the choice of the boundary conditions employed when defining the DH functions [see Eq.~(\ref{EQN:BC4})]. This is important as the gradient of the reconstructed functions on the free edge of the disk is imposed by the boundary condition choice. Finally, we demonstrate how the Fourier-Bessel PSD is invariant with regard to the curvature of the analysed surfaces and discuss the errors introduced by curvature projections.

\subsection{Fourier-Bessel functions and Hurst exponent}
\par
Persson et al.~\cite{Persson2005} showed that the decay of the amplitudes obtained with the Fourier PSD is equal to $-2(1+H)$ for an isotropic surface, where $H$ is the Hurst exponent. In \ref{APP:fract_PSD}, we derive an equivalent relationship between the decay of the PSD obtained from the Fourier-Bessel basis function and the Hurst coefficient. For this purpose, we use only the Bessel functions of the first kind and zeroth order ($m = 0$). This restriction is similar to the one employed in Fung (1967) \cite{Fung1967}, who generalised the traditional Fourier transform of 1D signals to 2D and 3D signals with a modified Wiener-Khintchine theorem dealing with Bessel functions of the first kind and zeroth order. The zeroth order is also appropriate for the isotropic self-affine surfaces that are analysed here due to the radial symmetry of these bases. For these surfaces, we show that the decay in amplitudes of the PSD spectrum of order $m = 0$ obtained with Fourier-Bessel bases is proportional to the Hurst coefficient as:
\begin{equation}
    D_{k, m = 0}^2 = \Big|\Big|F\Big(l(0)_k\Big)\Big|\Big|^2 \propto \lambda^{-2 (3/4+H)},
    \label{EQN:Main_PSD_vs_H}
\end{equation}
where $\lambda = l(0)_k$, which is a notation that will be used interchangeably. This is important because the direct relationship between decay and Hurst coefficient can be used to estimate the Hurst coefficient from the Fourier-Bessel PSD. There are two main advantages to using the Fourier-Bessel PSD rather than the standard Fourier PSD: (i) It can be applied to any layout of vertices, meaning the vertices do not need to lie on a predefined regular lattice (grid) and can be arranged arbitrarily. (ii) We can use this method to investigate the Hurst exponent of either nominally flat or curved isotropic self-affine surfaces without prior treatment of the surfaces (resampling of vertices and removing the tilt, see \cite{Jacobs2017}) or presumptions about the BC periodicity of the surface. Instead, we assume a Neumann boundary condition on the edge of the surface to allow for any arbitrary shift along the edge; we show later that this does not show influence when analysing the surfaces.
\par
\begin{figure*}
\centering
\includegraphics[width=1.0\textwidth]{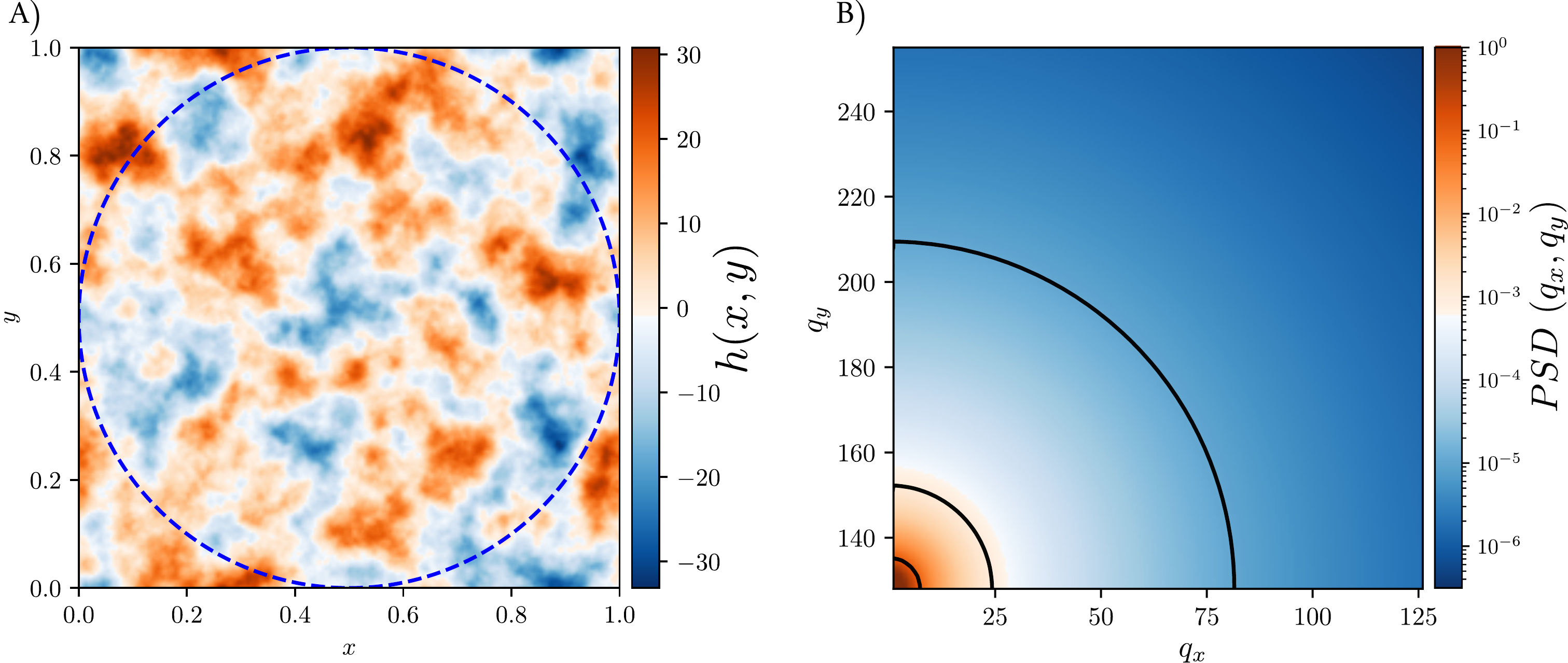}%
\caption{Randomly generated self-affine isotropic surfaces made using the Fourier filter method. A) A periodic fractal surface with Hurst exponent $H = 0.9$ and corresponding fractal dimension $D = 2.1$. The root mean square (RMS) of the height map $h(x,y)$ is $10$, where all the spatial units are measured in (length unit). B) Representation of the isotropic power law $C^{iso}(q) = q^{-2(1+H)}$ used for generating self-affine surfaces with $q_r = 0$, $q_l = 2^2$ and $q_s = 2^8 = 256$ with a seed number of $101$ (for reproducibility).
\label{FIG:PSD_SEFL_AFFINE}
}
\end{figure*}

This relationship between the Hurst exponent and the decay in the Fourier-Bessel PSD can be confirmed by testing fractal surfaces generated as described previously \cite{Jacobs2017, Persson2005}. The generation procedure, also described briefly in \ref{APP:generating_sefl_affine}, produces isotropic surfaces such as the one shown in Fig. \ref{FIG:PSD_SEFL_AFFINE}. This surface-generating code, like all other codes produced for this paper, is publicly available (see Section \ref{SEC:Repro}).

\subsection{Analysing flat and rough surfaces}
\par
In this sub-section, curvature-free surfaces are only considered. To determine the Fourier-Bessel PSD from a disk-shaped parameterisation space, it is natural to extract circular patches from surfaces. Sampling such patches must be done cautiously to prevent the loss of statistical information from the generated surfaces; for this reason, we sample, as shown in Figure \ref{FIG:PSD_SEFL_AFFINE}A, the largest circular patch that fits into the original surface.
\par
\begin{figure*}
\centering
\includegraphics[width=1.0\textwidth]{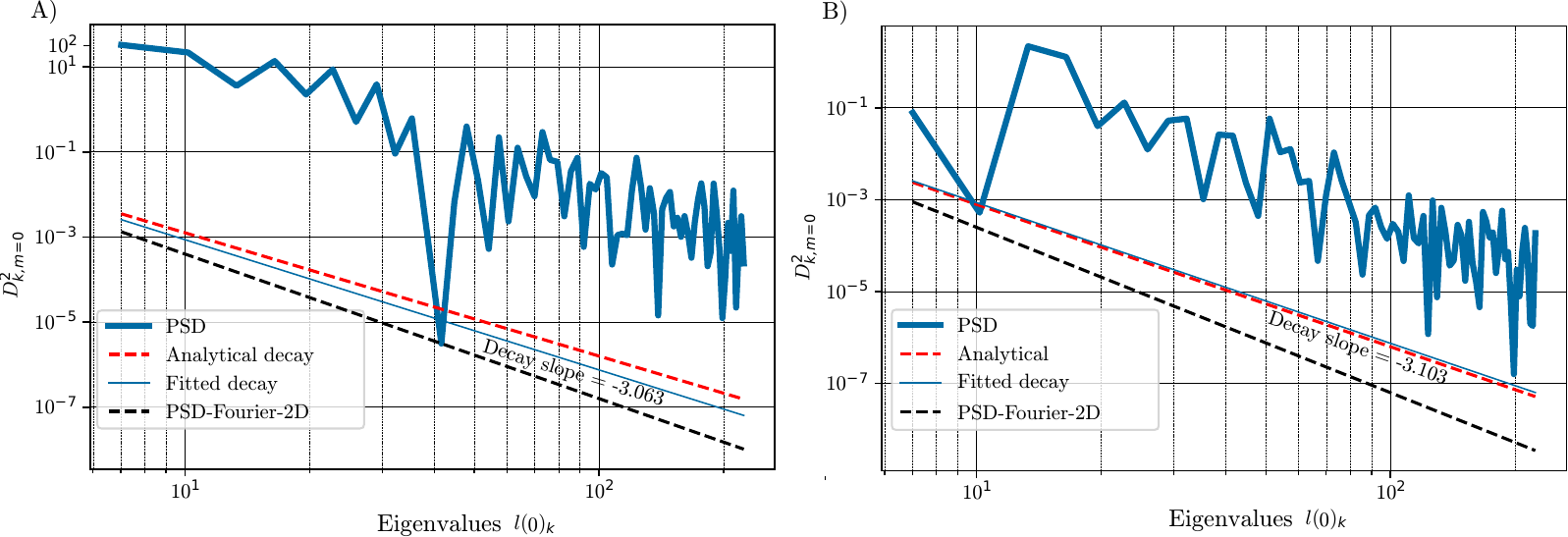}%
\caption{A) The power spectral density (PSD) computed using the Fourier-Bessel basis \big($D_{k \geq 2,m=0,z}^2$\big) with a seed number of $101$ (see Figure \ref{FIG:PSD_SEFL_AFFINE}) where the computed Hurst exponent of $H = 0.78$. B) The PSD \big($D_{k \geq 2,m=0,z}^2$\big) with a seed number of $51$ and a computed Hurst exponent of $H = 0.80$.
\label{FIG:PSD_07_and_08}
}
\end{figure*}

Figure \ref{FIG:PSD_07_and_08}A and B show the computed PSD for surfaces generated with Hurst exponents of $H=0.7$ and $0.8$, respectively. In this paper, the fitted slopes are obtained from fitting a power-law equation $D_{k,m=0}^2 = a \lambda^{b}$ to the Fourier-Bessel PSDs on a log-log scale; to do so, the least squares fit was used over the entire eigenvalue range that corresponds for $2 \leq k \leq 70$; the value of $l(0)_{70} > 200$ is close to the eigenvalue of $256$ of Fourier wherein we generated the surface. For an easier visual comparison, the best fits are indicated as shifted lines ($a$ was set to $1$) in order to separate from the PSD curves. In these plots, the line labeled \say{Analytical} has the perfect slope of $-2(3/4+H)$ [Eq.~(\ref{EQN:summary_PSD_vs_H})], where $H$ being the exponent used for the generation (i.e., $H=0.7$ and $0.8$ respectively). As the analysed surfaces are nominally flat, we considered only the z-component of the PSD as all the fractal information is embedded in that direction. The same plot shows that the fitted slopes produce good estimates of the decay and therefore of the Hurst exponent (error of $5.62\%$ and $0.10\%$ for surfaces in Fig.~\ref{FIG:PSD_07_and_08}A and B, respectively). The plots also show the analytical decay of the PSD from Fourier in the 2D case for comparison purposes.

\subsection{Error introduced by sampling a circular patch from a patch of a different shape}
\par
One source of error in estimating the Hurst exponent is the loss of data that results from sampling a circular patch from an original patch with a different shape. Non-circular patches that differ from a disk shape force the bases along the $x(\rho, \phi)$ and $y(\rho, \phi)$ to morph into a non-circular patch shape resulting in more artefacts when extracting the Fourier-Bessel amplitudes (sharp corners are expected to need a large number of modes to be captured from a disk parameterized disk). These artefacts do not affect the estimate of the Hurst exponent if the surface is flat because one may then consider only the PSD component computed along the $z$ direction. However, many real surfaces are curved surfaces. For these surfaces, the Hurst exponent is also embedded in the $x$ and $y$ coordinates of the surface. For such surfaces, sampling directly a circular patch rather than a non-circular patch leads to fewer artefacts when extracting the Fourier-Bessel amplitudes. Sampling a circular patch, however, leads to a loss of data.

\par
\begin{figure*}
\centering
\includegraphics[width=1.0\textwidth]{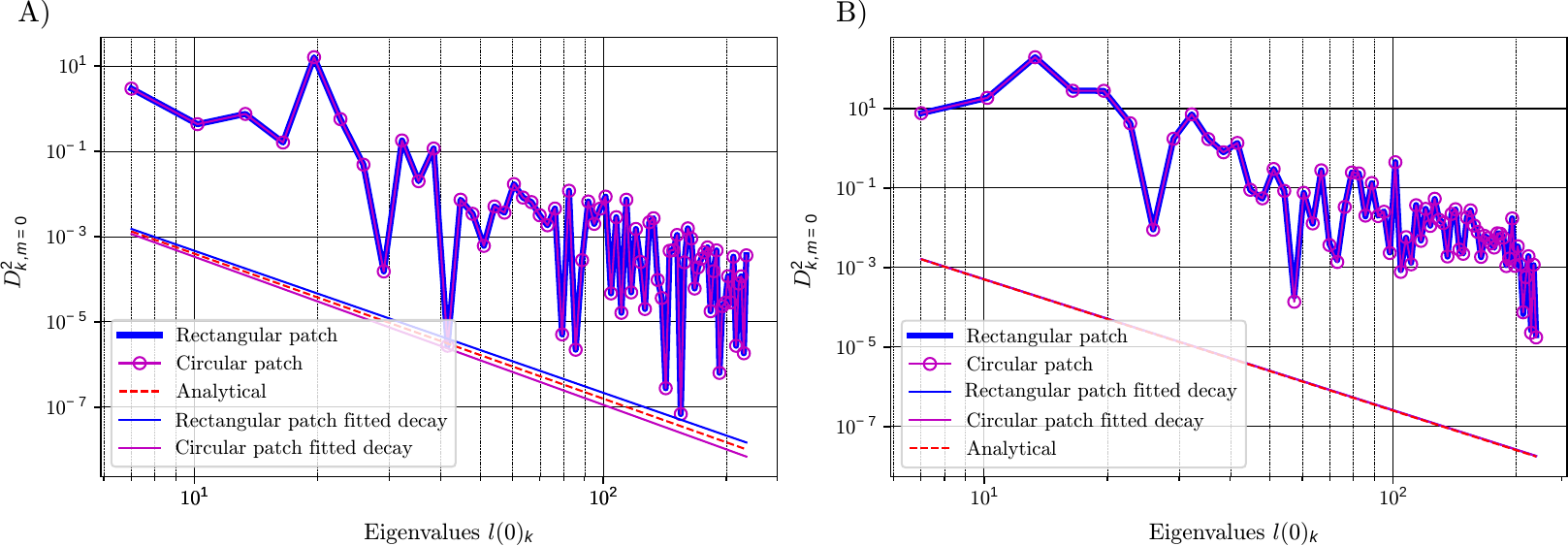}%
\caption{Comparison between square and circular patches, demonstrating that the chosen boundaries of the surfaces do not affect the PSD, assuming no loss of data in the sampled circular patch. A) The tested surface was generated with $H = 0.95$ and a seed number of $2345$. The obtained Hurst exponents were $0.92$ and $0.99$ for the square and circular patches (error of -3.15\% and +4.21\%). B) The tested surface was generated with $H = 0.9$ and a seed number of $7786$. The obtained Hurst exponents were $0.8967$ and $0.8965$ (error of -0.37\% and -0.38\%) for the square and circular patches.
\label{FIG:PSD_REC_VS_CIRC}
}
\end{figure*}

Here, we test how the PSD is sensitive to the loss of data that results from sampling a circular patch from a square. For this purpose, we generate square patches using the FFT-filter approach. First, we sampled the largest circular patch that fits into the generated square fractal surfaces (see Fig.~\ref{FIG:PSD_SEFL_AFFINE}A for illustration) and analysed it with the DHA. Second, we analysed the whole square patch and compared the computed PSDs in both cases along. We consider only the PSD component computed along the $z$ direction as we prior know that the surfaces are nominally-flat and free of curvature. In Fig. \ref{FIG:PSD_REC_VS_CIRC}, we show the extracted PSD for two types of surfaces ($H = 0.9$ and $0.95$) generated with various seeds for the random heights generator. The results show differences less than $0.20\%$ between the values obtained for the fitted decay and the Hurst exponents used for the surface generation. However, the Hurst exponents derived from the original square patch analyses were always slightly closer to the one used for the generation. This is expected as the circular patches do not consider all the surface data. Larger errors could be introduced while using even smaller circular patches, as more data are excluded when computing the PSD, especially of large wavelengths.
\par
In general, the resulting PSD curves, as shown in Figs.~\ref{FIG:PSD_07_and_08} and \ref{FIG:PSD_REC_VS_CIRC}, are noisy, likely due to the harmonic disagreement between the eigenvalues used to generate surfaces (via Fourier bases) versus the eigenvalues used to analyse them (Bessel-Fourier bases). Also, the discretization used for the DH parameterisation domain is not necessarily regular as was for the uniform grid employed for the Fourier analysis. In other words, the noise in the PSD here corresponds to the inter-harmonics, i.e. the frequency components in between the ones discretized with the discrete Fourier transform [see also the comments on Eq.~(\ref{EQN:PSD_H_Bessel}) about using different eigenvalues]. Such effects are particularly visible in artificial surfaces generated with FFT filters; in real rough surfaces, though not shown here, we have not observed these effects. Nevertheless, the results are in acceptable agreement with the expectations. Finally, the relatively low value (less than $10^{-5}$) as seen in the PSD in Fig.~\ref{FIG:PSD_07_and_08}A at about $l(0)_k \approx 40$ corresponds to uncorrelated interharmonics and is considered an outlier. Discarding such outliers can enhance the quality of the Hurst prediction (fitting).

\subsection{Influence of the applied Neumann boundary conditions on the PSD}
\par
\begin{figure*}
\centering
\includegraphics[width=1.0\textwidth]{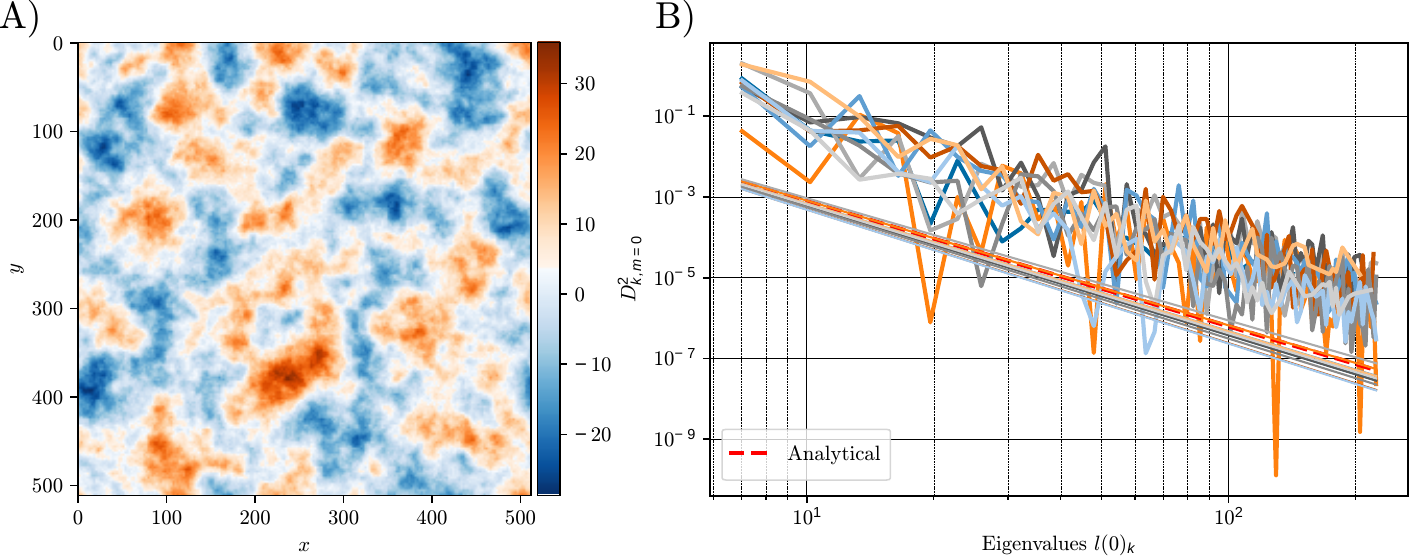}%
\caption{PSD comparison between multiple patches of the same size sampled over the same surface. A) The generated self-affine surface with $H=0.8$, $q_s = 2^9 = 512$ and a seed number of $807$. B) Comparison between the obtained PSD per randomly sampled patch and relative to the analytical prediction; fitting colours match the colours of the computed PSD.
\label{FIG:PSD_BC_samples}
}
\end{figure*}
\begin{figure*}
\centering
\includegraphics[width=0.6\textwidth]{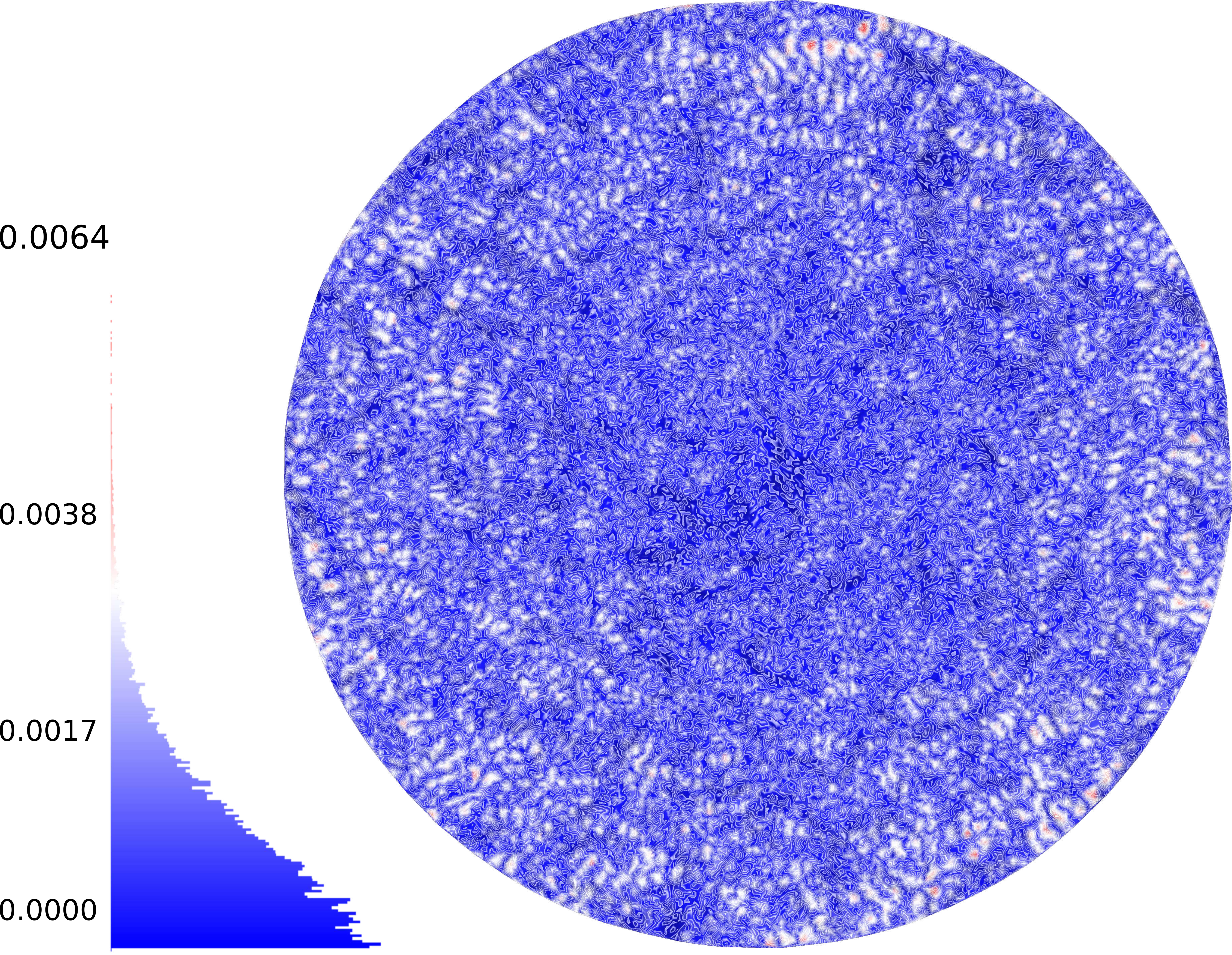}%
\caption{The surface reconstruction error map of the Hausdorff distance. The side histogram shows the normalised error distribution on the reconstructed surface with up to $k=71$ degrees in DHA with an edge length of $0.005$ (corresponds 
to $145,103$ reconstruction vertices) on a unit disk. This error is naturally occurring from fitting the smooth harmonic functions into a rough discontinuous surface. It can be noticed that the edge reconstruction error is consistent with the error distribution across the whole surface showing that the Neumann boundary conditions have no effect on the final results. 
\label{FIG:rec_map_error}
}
\end{figure*}

In Section \ref{SEC:FOURIER_BESSEL}, the Fourier-Bessel basis functions were computed by means of a specific Neumann boundary condition imposed on the outmost edge of the disk [see Eq.~(\ref{EQN:BC4})]. This allows arbitrary shifts along the edge of the disk but also constrains the gradients of the basis functions along the same edge, which can introduce artefacts if the surface profile is not compatible with such shapes. In the following we test whether the slope of the Fourier-Bessel PSD is sensitive to such a boundary condition by using several patches of the same size, placed at various locations on the self-affine fractal surface. To do so, a self-affine fractal surface four times as large as the one presented in Fig.~\ref{FIG:PSD_SEFL_AFFINE} was generated for a Hurst exponent $H = 0.8$ (see the large surface on Fig. \ref{FIG:PSD_BC_samples}A). Then, ten circular patches (of similar size to the above examples) are placed at random locations on the generated fractal surface such that the sampled disks are entirely within the boundary of the generated square surface. Then, the PSDs of these circular patches are computed using the DH analysis. Fig.~\ref{FIG:PSD_BC_samples}B shows the obtained PSDs and their corresponding slopes. As already shown in Fig.~\ref{FIG:PSD_REC_VS_CIRC}B, the circular patches lead to a very good estimate of the Hurst coefficient of this surface with an average slope of $-3.193 \pm 0.085$ (mean $\pm$ standard deviation) and the expected decay was $-3.1$, corresponding to $H=0.8$ for the Hurst exponent. If the Neumann BC of the Fourier-Bessel basis function would have a large influence on the results, a significant error scatter in the Hurst coefficient estimates would be obtained since the mismatch at the disk edge varies with the position of the sampling disk. We can therefore argue that the Neumann BCs have a non-significant influence over the Hurst extraction. In fact, the Neumann BC imposed at the edge of the parameterisation domain $\rho(\phi) = 1$ is only affecting the endpoint, and the slope in the vicinity of that edge is arbitrarily changing with the bases and the corresponding weights $q_m^k$ (see Fig.~\ref{FIG:BESSEL} and the discussion in \ref{APP:BCs}). If the Neumann BC of the Fourier-Bessel basis function has a significant effect on the computed PSD, the reconstructed surface should have larger errors at the edge than in the centre. Figure \ref{FIG:rec_map_error} shows the error of such reconstructed surfaces, which shows that the error is not significantly increased towards the edge of the surface. The shown error map has been normalized by the diagonal length of the smallest bounding box of the rough surface. The reconstruction error is the result of many factors such as the used mapping and regularisation steps. However, the choice of the basis functions can affect the error as was explained by \cite{Verrall:98}.

\subsection{Curvature normalisation of rough surfaces}
\par
In this section, we show that the computed PSD is invariant to surface curvature as fully described in \ref{APP:analytic_cap}. This is essential for studying arbitrarily curved rough surfaces. In general, when studying a sub-part of a surface (usually when using any microscopy acquisition), the curvature may be hardly visible even though it would pollute the PSD. Computing a PSD that is invariant to the curvature is essential to enable multi-scale surface analysis.
\par
\begin{figure*}
\centering
\includegraphics[width=0.8\textwidth]{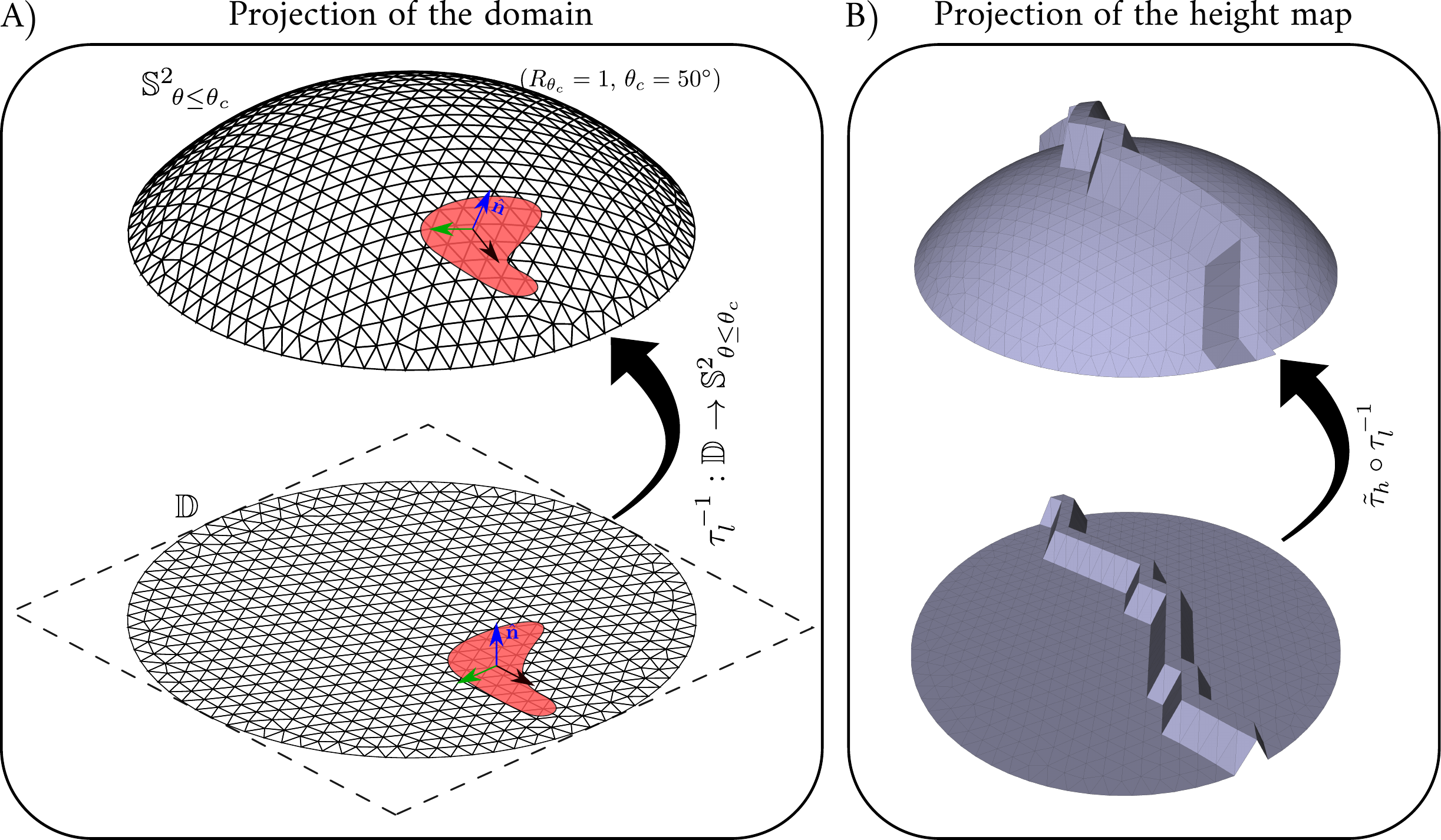}%
\caption{Schematic of the curvature projection. A) The inverse Lambert projection $\tau^{-1}_{l}$ between smooth surfaces. B) The radial projection $\tilde{\tau}_h$ for the height maps $h(x,y)$ such that the asperities remain normal to the surface, as illustrated with a random stroke of constant height.
\label{FIG:SCHMTC_CURV_PROJ}
}
\end{figure*}

The procedure employed to generate curved and self-affine rough surfaces is now presented. The level of curvature was controlled by projecting nominally flat surfaces onto spherical caps (see \ref{APP:generating_sefl_affine} for more details). Spherical caps were used as they have a uniform curvature radius $R_{\theta_c}$, where $\theta_c$ is a chosen half-angle controlling the desired curvature. The Lambert azimuthal inverse projection was used to obtain a mapping from the sampled circular patch towards the spherical patch. Then, the rough height map is radially projected onto the spherical cap, as shown in the schematic in Fig.~\ref{FIG:SCHMTC_CURV_PROJ}. Such a radial projection implies that asperities are normal to the curved surface, ensuring that the curved height maps are not restricted to only having $z$ components. However, this projection leads to distortions of the angles in the meshing triangles. The effects of such distortions will be discussed next.
\par
\begin{figure*}
\centering
\includegraphics[width=1.0\textwidth]{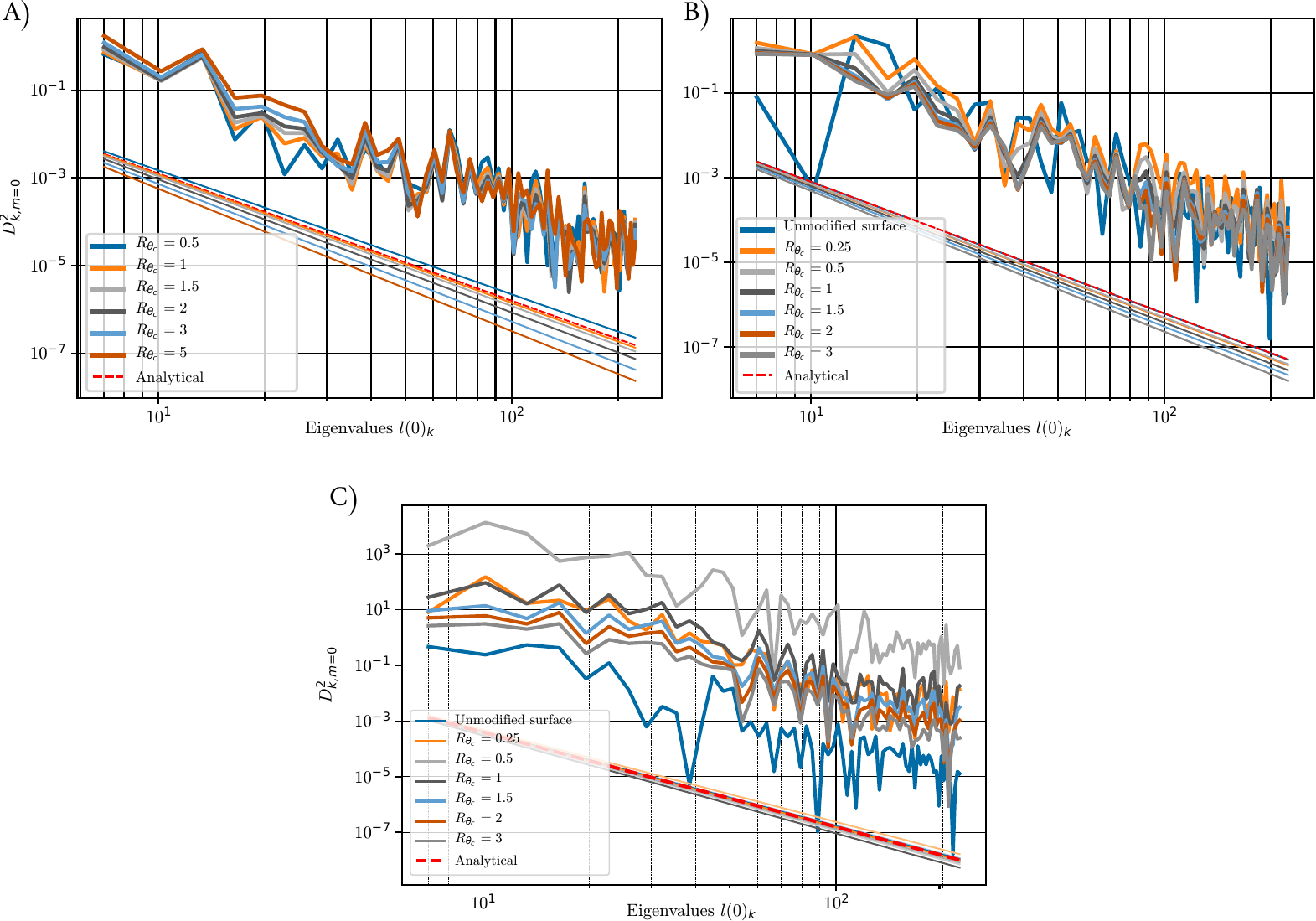}%
\caption{Effects of curvature on rough spherical patches with various radii $R_{\theta_c}$ ($\theta_c = 10^{\circ}$) and $H\in \{0.7, 0.8, 0.95\}$ with seed numbers of $\{101, 51, 9700\}$, respectively. A) The normalised PSD for the curvature \big($D_{k \geq 2,m=0}^2$\big) as explained in Eqs.~(\ref{EQN:SHAPE_DESC}) and (\ref{EQN:NORM_CAP_DESC}) for the generated surface in Fig.~\ref{FIG:PSD_SEFL_AFFINE}, where the fitted PSD considers $2\leq k \leq 70$. The fitted slope for various spherical patches, which were $\{-2.827, -2.925, -2.961, -3.034, -3.138, -3.247\}$ for $R_{\theta_c} = 0.5, 1, \dots, 5$, respectively; the analytical slope for $H=0.7$ is $-2.9$. B) Surface generated with $H=0.8$ with $R_{\theta_{c}} = \{\infty, 0.25, 0.5, 1, 1.5, 2, 3\}$. Respectively, fitted decay slopes here were $\{-3.103, -3.159, -3.149, -3.208, -3.257, -3.316, -3.362\}$ considering that the analytical slope for $H=0.8$ is $-3.1$. Notice that we here denote the unmodified surfaces (flat configuration) by $R_{\theta_{c}} \to \infty$. C) Surface generated with $H=0.95$ and $R_{\theta_{c}} = \{\infty, 0.25, 0.5, 1, 1.5, 2, 3\}$. The fitted decay slopes here were $\{-3.373, -3.311, -3.467, -3.389, -3.406, -3.438, -3.517\}$ for $R_{\theta_c} = \{\infty, 0.25, \dots, 3\}$, respectively, considering that the analytical slope for $H=0.95$ is $-3.4$. Note that the fitting colours match the colours of the computed PSD.
\label{FIG:CURV_PSD_VS_NORM}
}
\end{figure*}

We applied this approach to circular patches, of unit disks, that were generated with input Hurst exponents $H=0.7$, $0.8$ and $0.95$. We projected the surfaces onto several spherical caps with radii ranging from $0.25$ up to $5$ (unit length). The size of the cap is controlled via the half-angle $\theta_c$, while the radius of the cap controls directly the constant Gaussian curvature of the cap. We then performed a PSD analysis with Fourier-Bessel basis functions and then fitted the decaying slope to the moments ${D}_{k}^2$; using the resultant shape descriptors along $x$, $y$ and $z$ and not only the $z$ one as was done for nominally flat surfaces. Figure \ref{FIG:CURV_PSD_VS_NORM} shows the obtained Hurst exponents obtained from spherical caps of various radii. It can be seen that the mean decay slopes obtained from curved caps were $-3.421\pm0.064$ for $H=0.95$, $-3.242\pm0.078$ for $H=0.8$ and $-3.022\pm0.138$ for $H=0.7$. Hurst exponents from the mean decay slopes correspond to mean Hurst estimates of $0.96$, $0.87$ and $0.76$. The following observations can be made:
\begin{itemize}
    \item In Fig.~\ref{FIG:CURV_PSD_VS_NORM}, we show the angular distortion due to the projection of the flat disks into caps. Good estimates of fractal surfaces are normally obtained when the distortion due to mapping is minimal. However, as the used projection is not isometric, the projection onto the curved surface does seem to introduce spurious noise, most likely due to the area-preserving algorithm which does not preserve angles or distances (see Section~\ref{APP:roughness_proj} for further details).
    
    \item In general, Hurst exponents of smoother curved surfaces (larger Hurst exponents) tend to be better estimated than for rougher curved surfaces. This is again attributed to the distortion introduced by the projection (see Fig.~\ref{FIG:angle_distortion}). Overall, this provides numerical evidence that the method presented in this paper is capable of computing curvature-free PSD of surfaces assuming there is a method to generate fractals on curved surfaces without projection artefacts or loss of data.
\end{itemize}
\par
It also should be noted that the projection approach distorts the fit for the few first eigenvalues. Generally speaking, the larger the cap radius (whence $R_{\theta_c} \geq 1$), the larger the distortion in our PSD and the corresponding computed decay (see the brief explanation in \ref{APP:roughness_proj} and Fig.~\ref{FIG:angle_distortion} about how we measured and defined the angular distortion). Regardless, the purpose of this paper is not to handle and correct the projected distortion on such caps, but to show the ability of the newly proposed method to be insensitive to curvature. 

%% file: morph_results.tex
\par
\begin{figure*}
\centering
\includegraphics[width=1.0\textwidth]{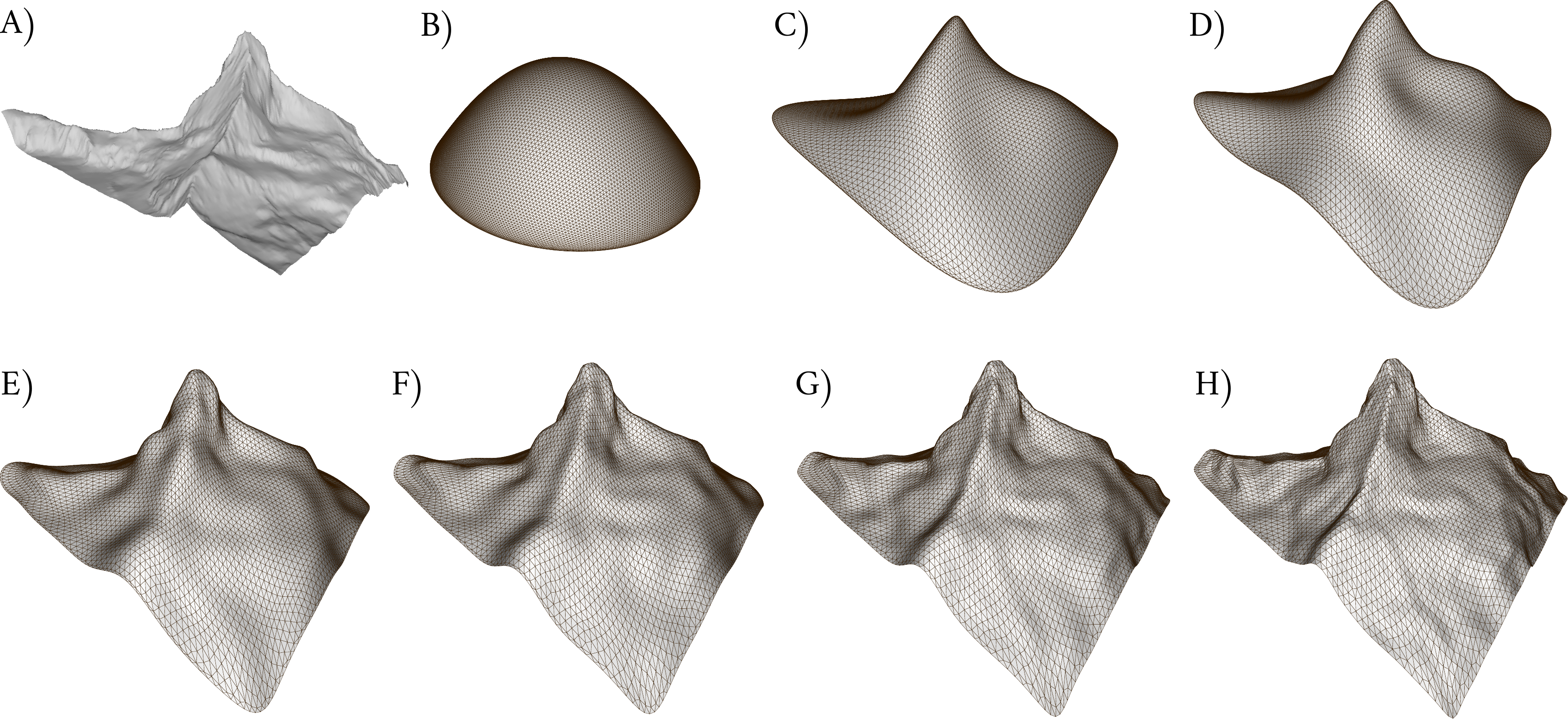}%
\caption{DH reconstruction of the Matterhorn mountain in Zermatt, Switzerland. A) The northeast view of the Matterhorn \cite{REF:9, REF:10}. The mountain was expanded up to $k = 50$ and reconstructed with an edge-length of $0.025$ ($5,809$ reconstruction vertices) on a unit disk. The reconstruction results at different degrees: B) $k=1$ results in an FDEC with a size of $2a = 601.77, 2b= 512.18 \text{ and } c = 176.97$; C) $k = 3$; D) $k = 5$; E) $k = 15$; F) $k = 20$; G) $k = 30$ and H) $k = 50$. Note $k=5$ means that all degrees up to 5 are included.
\label{FIG:Matterhorn_rec}
}
\end{figure*}

Descending from the spherical harmonics and elliptical Fourier descriptors, the DH method we propose here is useful for numerous applications beyond the physics of rough surfaces, such as for decomposing and reconstructing surfaces for statistical studies, including medical imaging. To demonstrate such applications, we herein propose the Matterhorn mountain in Zermatt, Switzerland as a benchmark for genus-0 single-edge open surfaces. Figure \ref{FIG:Matterhorn_rec}A shows the input surface consisting of $93,635$ vertices, and the following images (B-H) show reconstructions at various degrees. In general, the higher the decomposition degree, the smaller the wavelength, and thus, the higher the level of spatial resolution in the reconstruction.
\par
Figure \ref{FIG:Matterhorn_rec}B shows the reconstruction of the Matterhorn from the FDEC ($k = 1$) up to the inclusion of the first $50$ degrees ($k = 50$). To generate unit disks with uniform point distributions, for the reconstruction problem, we used the method explained in Shaqfa et al. (2021b) \cite{shaqfa2021b}. Setting an element edge size of $0.025$ on the reconstruction domain of the unit disk provided $5,809$ uniformly distributed reconstruction vertices. The error between the input and reconstructed surface was measured by the RMSE of the Hausdorff distance, normalised by the length of the diagonal distance of the bounding box. For the reconstructed surface with $k = 50$, the error was $0.000703$, which is very small and means that the reconstruction matches the input surface. However, an edge size of $0.0075$ of a unit disk (with $64,492$ reconstruction vertices) produces an RMSE of $0.000628$ (not shown on the figure). The RMSE of the 2-norm distance between the input and output vertices for edge lengths of $0.025$ and $0.0075$ were $0.402921$ and $0.360014$, respectively. Such results demonstrate the robustness and suitability of the proposed approach for genus-0 single-edge open surfaces. More detailed comparisons between the SCHA \cite{shaqfa2021b} and DHA methods as well as a brief commentary can be found in \ref{APP:SCHA_vs_DHA_res} using a 3D face as a benchmark.

%% file: conclusions.tex
\par
To address the lack of methods for analysing the roughness and morphology of nominally flat and open surfaces, in this paper, we used disk harmonics (DH) with the Fourier-Bessel basis function to decompose genus-0 single-edge disk topologies. We first combined the density-equalising and quasi-conformal maps to compute a bijective area-preserving parameterisation of open surfaces onto a unit disk, which allows us to use DH for these surfaces. The proposed DH approach is not only useful for analysing and reconstructing surfaces, such as in biomedical imaging, but we also demonstrated its ability to measure the Hurst exponent (fractal dimension) of isotropic self-affine fractal surfaces. The latter was reached by exploiting an expansion of the Wiener-Khintchine theorem that links the autocorrelation functions to the PSD computed by the DH basis of the first kind and zeroth order. The obtained relationship between the PSD and the DH basis (of $m=0$) were confirmed by generating artificial fractal surfaces using the Fourier-Filter method. 
To study the influence of the curvature on the measured PSD, we generated rough fractal surfaces with prescribed curvatures using the Lambert azimuthal inverse projection to map sampled circular patches onto spherical caps. These rough caps were then used to test the independence of the Fourier-Bessel PSD from the principal curvatures of the surface. Overall, our results correspond well with the theoretical expectations presented in this paper. This approach should overcome many of the drawbacks that can be found in the standard Fourier-filter approach, which implies the use of periodic boundary conditions as well as data interpolated on a regular rectangular lattice. Accounting for curvature is a new feature that will improve the modelling of more complicated problems in the future, such as rough contacts between real scanned surfaces and the mechanics of complicated interfaces. We also envisage that this paper will assist in generating more realistic rough surfaces, besides many other aspects beyond the purpose of this paper.

%% file: reproducibility.tex
The codes and the used surfaces are made available in this paper and can be found on GitHub: \href{https://github.com/eesd-epfl/disk-harmonics}{Disk harmonics}.

%% file: appendix.tex

\subsection{Estimating the spatial wavelengths associated with Fourier-Bessel functions}
\label{APP:wavelengths}
\input{app_A_wavelenghts}

\subsection{Curvature normalisation of the parametric form of an ellipsoidal cap}
\label{APP:analytic_cap}
\input{app_A_curvature_normalisation}

\subsection{The size of the first degree ellipsoidal cap (FDEC)}
\label{APP:FDEC_size}
\input{app_A_FDEC_size}

\subsection{Hurst exponent of an isotropic self-affine surface with Fourier-Bessel functions}
\label{APP:fract_PSD}
\input{app_A_fract_PSD_analytical}

\subsection{Isotropic power-law for self-affine scaling surfaces}
\label{APP:generating_sefl_affine}
\input{app_A_generating_self_affines}

\subsection{Generating curved fractal surfaces over spherical caps}
\label{APP:roughness_proj}
\input{app_A_roughness_projection.tex}

\section{Supplemental results}

\subsection{A commentary on Bessel function boundary conditions}
\label{APP:BCs}
\input{app_B_boundaryconditions}

\subsection{Eigenvalues of Bessel functions of the first kind}
\label{APP:Eiegn}
\par
\begin{figure*}
\centering
\includegraphics[width=0.7\textwidth]{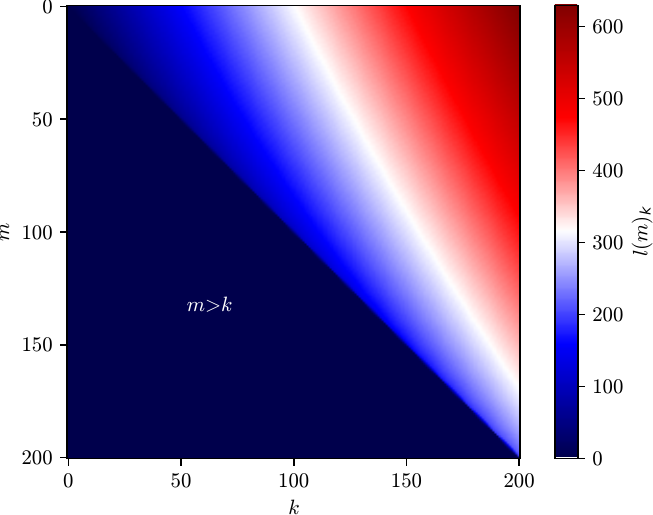}%
\caption{Eigenvalues for Bessel functions for $k\in\{0, 1, \dots, 200\}$. Notice when $m>k$, $l(m)_k=0$.
\label{FIG:EIGEN_K_200}
}
\end{figure*}

Many methods can be used for computing the eigenvalues (zeroes) of Bessel functions with the proper boundary conditions. In this paper, we used the same approach we used in the SCHA method by Shaqfa et al. (2021) \cite{shaqfa2021b}. Figure \ref{FIG:EIGEN_K_200} shows the roots (zeros) identified for the first $200$ degrees.


\subsection{Comparison between DHA and SCHA for reconstructing open surfaces}
\label{APP:SCHA_vs_DHA_res}
\input{app_B_comparison_DHA_vs_SCHA}

%% file: app_A_wavelenghts.tex
\par
\begin{figure*}[!ht]
    \centering
    \includegraphics[width=1.0\textwidth]{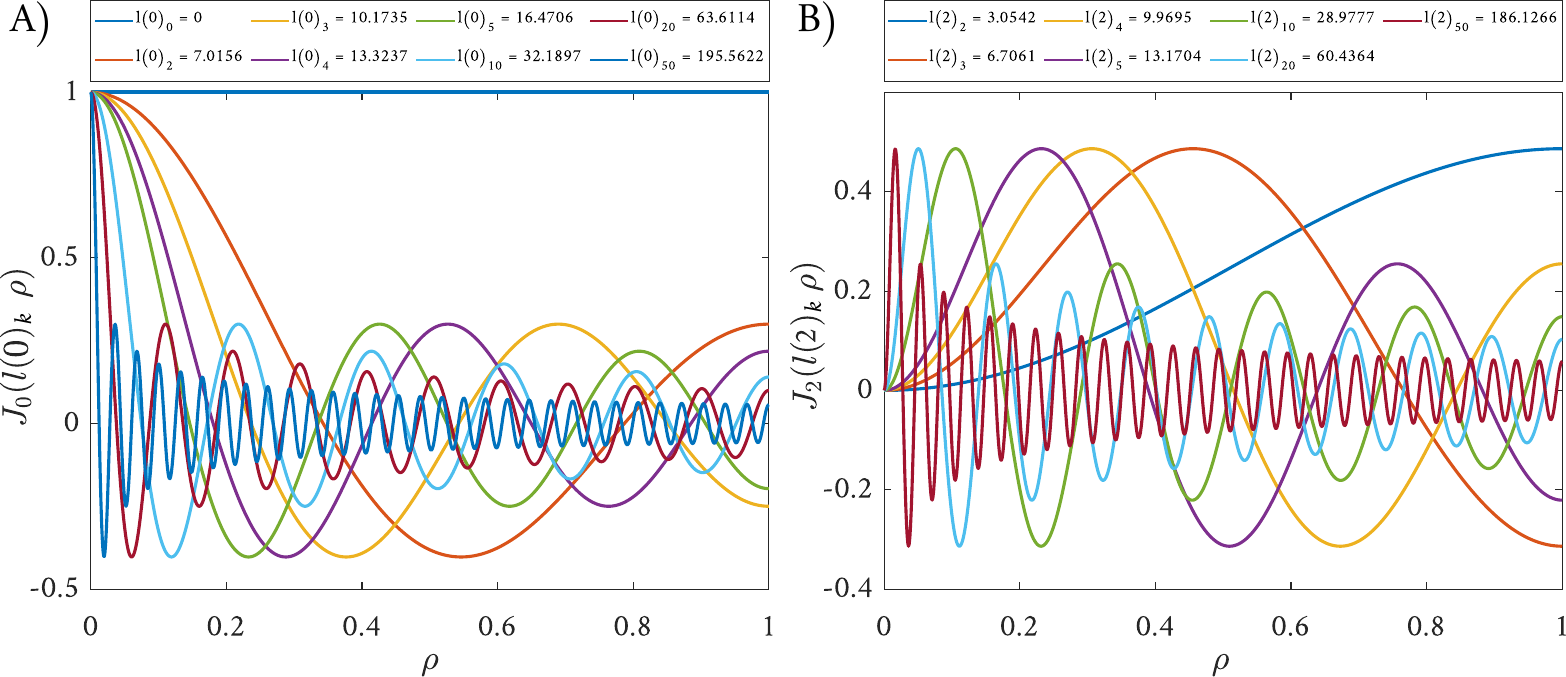}
    \caption{The unnormalized Bessel solutions for $\rho \in [0, 1]$, which change signs $k-m$ times within $\rho \in [0, 1]$. Notice the influence of the Neumann boundary conditions near $\rho \leq 1$ as ${\partial J_m (l(m)_k \cdot 1)}/{\partial \rho} = 0$ (cf. with Fig.~\ref{FIG:BESSEL_Dirichlet_BC}, with Dirichlet BC).}
    \label{FIG:BESSEL}
\end{figure*}

The frequency with which Eq.~(\ref{EQN:LAP_BESSEL}) changes its sign within $\rho \in[0, 1]$ depends on $l(m)_k$: the function changes sign $k-m$ times within $\rho \in [0, 1]$, as can be seen on Fig.~\ref{FIG:BESSEL}in the case where $m=0$ and $m=2$. This is closely related to the spatial frequency inherent to the basis functions. 
\par
However, to estimate the wavelengths associated with any given surface, termed the spatial data, and link it to the corresponding analysis degrees, it is natural to assume that the order $m$ is very small: $m \to 0$ and $k \gg m$ as $k \to \infty$. These assumptions are common in this field; see our previous work \cite{shaqfa2021b} for details about using such assumptions in spherical harmonics and spherical cap harmonics. Here, we derive two expressions for the wavelengths from both parts of the harmonic solution: (i) from Bessel functions by considering radial wave variation and (ii) Fourier functions with the angular wave variation.
\par
As we assume that $k \to \infty$, we need to find the asymptotic behaviour of Bessel functions for large terms. Such behaviour can be found in traditional mathematics books \cite{REF:1}, where the Bessel function takes the form:
\begin{equation}
J_{m}(u) \approx {\sqrt {\frac {2}{\pi u}}}\sin \left(u-{\frac {m \pi }{2}}+{\frac {\pi }{4}}\right) + \mathcal{O}(u^{-3/2}) + \dots
\label{EQN:Bessel_ASYM}
\end{equation}
From Eq.~(\ref{EQN:Bessel_ASYM}), we can see that the Bessel function asymptotically behaves like a regular $\sin x$ wave, and its roots approach the zeros of the trigonometric function with a slight phase shift. We can also conclude that for large degree expansions, the spaces between eigenvalues monotonically approach $\pi$, such that $l(m)_{k+1}-l(m)_{k} \approx \pi$ and can be estimated by \cite{REF:1}:
\begin{equation}
    l(m)_k \approx k \pi + \frac{m\pi}{2} - \frac{\pi}{4}.
\end{equation}
Additionally, as $m\to0$, we can estimate $l(0)_k$ to be:
\begin{equation}
    \hat{l}(0)_k = k \pi - \frac{\pi}{4}.
\end{equation}
\par
Now, the wavelength of the radial functions that extend from the origin of the disk can be written as $\omega_{k}^{\rho} = 2 \pi r/\hat{l}(m)_k$; this term needs to be computed along the three principal directions $x$, $y$ and $z$. If we agree to consider the resultant of these terms and to consider the FDEC dimensions as the accepted domain size, we can write the following equations:
\begin{equation*}
    \omega_k^\rho \pi \Big(k - \frac{1}{4}\Big) \approx 2\pi\sqrt{a^2+b^2+c^2},
\end{equation*}
\begin{equation}
    \omega_k^\rho \approx \frac{2\sqrt{a^2+b^2+c^2}}{\Big(k - \frac{1}{4}\Big)}, \forall k \geq 1~\&~k \in \mathbb{Z}^{+}.
\label{EQN:RADIAL_WAVE_LENTH}
\end{equation}
Geometrically, Eq.~(\ref{EQN:RADIAL_WAVE_LENTH}) assumes equidistant oscillations and then takes the average variation along a line that passes through the origin of the disk connecting any two arbitrary points that lay on the edge of the disk (antipodal points).
\par
Now, consider the angular portion of the basis functions that has a minimal wavelength of $\omega_{min}^{\phi} = \xi/m$, where $\xi$ is the length of the perimeter of the outermost closed-edge of the FDEC, which can be approximated by $\xi = 2\pi \sqrt{({a^2+b^2})/{2}}$ with $a$ and $b$ representing the half-lengths defined on FDEC. If we assume that the smallest angular wavelength covered on the surface occurs when $k = m$, then the angular wavelength can be written as:
\begin{equation}
    \label{EQN:WAVELENGTH}
    \omega_k^\phi \approx \frac{2\pi}{k}\sqrt{\frac{a^2 + b^2}{2}}, \forall k \geq 1.
\end{equation}
Equation~(\ref{EQN:WAVELENGTH}) is equivalent to the term derived in Shaqfa et al. (2021) \cite{shaqfa2021b} that geometrically measures the variations (propagation) occurring only along the perimeter of FDEC. In this paper, we use Eq.~(\ref{EQN:WAVELENGTH}) as the wavelength $\omega_k$ and use the terms wavelength $\omega_k$ and angular wavelength $\omega_k^{\phi}$ interchangeably, as this angular estimation is the common factor in spherical harmonics, hemispherical harmonics, spherical cap harmonics and disk harmonics analyses. Additionally, the angular wavelengths are uniformly distributed over the length of the in-plane ellipse, though the radial wavelengths are not uniformly distributed along the radial lines of the disks or caps, as the radial variations accumulate at the edge of the disk for large values of $m$.

%% file: app_A_curvature_normalisation.tex
\par
The power spectrum is impacted linearly by the global curvature of the surface. In order to demonstrate this, the first order disk-harmonics decomposition needs to be mapped to an  ellipsoidal cap. 
With the convertion between Cartesian $(x,y,z)$ and cylindrical $(\rho, \phi, z)$ coordinates, i.e. $x = \rho \cos\phi$, $y = \rho \sin\phi$ and $z = z$, an ellipsoid centred over the origin axis will obey the following characteristic equation:
\begin{equation}
    \frac{x^2}{a^2} + \frac{y^2}{b^2} + \frac{z^2}{c^2} = 1,
    \label{EQN:ELLIPS}
\end{equation}
with $a, b, c$ the ellipsoidal characteristic lengths. The polar parameterisation of the unit disk
\begin{equation}
    x = a \rho \cos\phi;~y = b \rho \sin\phi
    \label{EQN:coord_conv}
\end{equation}
will be used, where the domain of this parametric space is $\rho \in [0,1]$ and $\phi \in [0, 2\pi]$.
Therefore, plugging Eq.~(\ref{EQN:coord_conv}) into Eq.~(\ref{EQN:ELLIPS}) gives:
\begin{eqnarray}
  \frac{a^2 \rho^2 \cos^2\phi}{a^2} + \frac{b^2 \rho^2 \sin^2\phi}{b^2} + \frac{z^2}{c^2} = 1 \nonumber \\
  \rho^2 + \frac{z^2}{c^2} = 1.
  \label{EQN:ELI2}
\end{eqnarray}
Then, the ellipsoid is parameterised with the equation:
\begin{equation*}
    z = c \sqrt{1 - \rho^2},
\end{equation*}
The final parametric form of the ellipsoidal cap becomes:
\begin{equation}
f(\rho, \phi) = 
\left\{
\begin{array}{l}
      x(\rho, \phi) = a \rho \cos\phi\\
      y(\rho, \phi) = b \rho \sin\phi\\
      z(\rho, \phi) = c \sqrt{1 - \rho^2}\\
\end{array}
\right.,
\label{EQN:ELLIPSE_PARAM}
\end{equation}
or in the general matrix form with an arbitrary rotation $R$:
\begin{equation}
\begin{pmatrix}
x(\rho, \phi) \\
y(\rho, \phi) \\
z(\rho, \phi) 
\end{pmatrix} = R \cdot \begin{pmatrix}
                          a \rho \cos\phi \\
                          b \rho \sin\phi \\
                          c \sqrt{1 - \rho^2}
                        \end{pmatrix}.
\label{EQN:ELLIPSE_PARAM_MAT}
\end{equation}
For the semi-axes of ellipsoids whose directions coincide with the Cartesian coordinate axes, $R$ is the identity matrix, which will be assumed for a moment. With this parametric form, giving the ellipsoid cap equation in terms of $\phi$ and $\rho$, we can expand the ellipsoidal cap and therefore compute the coefficients directly from the DH expansion. For instance along the x-axis any coefficient $q_{m,x}^1$ can be computed as:
\begin{equation}
  \frac{q_{m,x}^1}{a} = \int_0^{2\pi} \int_0^1 \rho^2 \cos\phi ~ D_m^1 ~ d \rho~ d \phi.
\end{equation}
Therefore the ratio $q_{m,x}^1/a$ happens to be a constant as it is equal to a definite integral. The same can be made with $q_{m,y}^1/b$ and
$q_{m,z}^1/c$. As a matter of fact, a similar relation holds for any vector coefficient $\vec{q}_{m}^k$, which suggests a linear relation between curvature and coefficients. This fact can be exploited to normalise, in each direction, a curvature-free PSD. The next appendix section deals with the general case where the principal axes must be extracted to compensate for the arbitrary rotation applied to the ellipsoid.

%% file: app_A_FDEC_size.tex
The general idea of fitting an ellipsoidal equation from a parameterised decomposition, draws from traditional principal component analysis (PCA), which solves for the principal axes and the associated lengths $a, b, c$ and the associated directions. To this end, we start with a decomposition limited to the first order $k=1$ only:
\begin{equation}
  \vec{f}^1(\rho, \phi) = \vec{q}_{-1}^1 D_{-1}^1(\rho, \phi) + \vec{q}_0^1 D_0^1(\rho, \phi) + \vec{q}_1^1 D_1^1(\rho, \phi).
  \label{EQN:besselfourier-approx}
\end{equation}
Here we want to assume that in the region $\rho \to 0$ the surface has a constant curvature that represents the curvature of the whole cap. Therefore, the Bessel functions can be expressed as a truncated series, for the first degree ($m \in \{ -1, 0, 1\}$):
\begin{equation}
    \label{EQN:bessel-approx1}
    J_0(\alpha \rho) = 1 - \left(\frac{\alpha \rho}{2}\right)^2 + \mathcal{O}\Big(\rho^4\Big) \qquad \forall \alpha \in \mathbb{R},
\end{equation}
\begin{equation}
    \label{EQN:bessel-approx2}
    J_1(\alpha \rho) = \frac{\alpha \rho}{2} + \mathcal{O}\Big(\rho^3\Big) \qquad \forall \alpha \in \mathbb{R}.
\end{equation}
From the relationship between the polar and Cartesian coordinates, we can obtain the approximation of the FDEC when close to the origin with $D_m^1$ approximated by using Eqs.~(\ref{EQN:bessel-approx1}) and (\ref{EQN:bessel-approx2}):
\begin{align}
  \hat{\vec f}(\rho, \phi) &= 
                             \left(
                             \begin{matrix}
                               - \frac{1}{2} N_{1}^1 l(1)_1 \rho e^{-i\phi} \\
                               N_{0}^1 (1+l(0)_1^2\rho^2) \\
                               \frac{1}{2} N_{1}^1 l(1)_1 \rho e^{i\phi}
                             \end{matrix}
                             \right)^T \cdot \left(
                             \begin{matrix}
                               \vec{q}_{-1}^1\\
                               \vec{q}_{0}^1\\
                               \vec{q}_{1}^1     
                             \end{matrix}
                             \right) + \mathcal{O}\Big(\rho^3\Big) \nonumber \\
                           &= \frac{N_{1}^1 l(1)_1 \rho}{2} (\vec{q}_1^1 e^{i\phi} - \vec{q}_{-1}^1 e^{-i\phi}) + N_{0}^1 \vec{q}_{0}^1 \Big( 1 + l(0)_1^2\rho^2\Big) + \mathcal{O}\Big(\rho^3\Big) \nonumber \\
                           &= \frac{N_{1}^1 l(1)_1 \rho}{2} \left[ (\vec{q}_1^1 - \vec{q}_{-1}^1) \cos(\phi) + i (\vec{q}_1^1 + \vec{q}_{-1}^1) \sin(\phi)\right]   + N_{0}^1 \vec{q}_{0}^1 \Big( 1 + l(0)_1^2\rho^2\Big)  + \mathcal{O}\Big(\rho^3\Big) \nonumber  \\
                           &=
                             \underbrace{ \frac{1}{2}  l(1)_1 N_{1}^1  \left(
                             \begin{matrix}
                                \vec{q}_1^1 - \vec{q}_{-1}^1\\
                                i( \vec{q}_1^1 + \vec{q}_{-1}^1)     
                             \end{matrix}
                             \right)^T}_{A^T}
                             \cdot 
                             \left(
                             \begin{matrix}
                               \rho \cos \phi \\
                               \rho \sin \phi
                             \end{matrix}
                             \right)
                             + N_{0}^1 \vec{q}_{0}^1 \Big( 1 + l(0)_1^2\rho^2\Big) + \mathcal{O}\Big(\rho^3\Big) \nonumber \\
                           &=  A^T \cdot \left(
                             \begin{matrix}
                               \hat{x}(\rho, \phi) \\
                               \hat{y}(\rho, \phi) 
                             \end{matrix}
                             \right) +  N_{0}^1 \vec{q}_{0}^1 \Big( 1 + l(0)_1^2\rho^2\Big) + \mathcal{O}\Big(\rho^3\Big), 
\label{eqn:first-order-fdec-form}
\end{align}
where
the symmetry $N_1^1 = N_{-1}^1$ and standard trigonometry
has been used to obtain the above expression. In this expression
$\hat{x}, \hat{y}$ are the Cartesian coordinates over the parameterised unit disk.
It must be noted that $A$ is a matrix since it is here expressed as a vector with components involving of the vectors $\vec{q}_1^1$ and $\vec{q}_{-1}^1$. Equation (\ref{eqn:first-order-fdec-form}) has a constant shift with the term $f(0, \phi) = N_{0}^1 \vec{q}_{0}^1$.
Without loss of generality for the study of curvature, it is possible to discard such a translation term, leading to the first order approximation:
\begin{equation}
  \hat{\vec f}(\rho, \phi) =  A^T \cdot \left(
                             \begin{matrix}
                               \hat{x} \\
                               \hat{y} 
                             \end{matrix}
                             \right) +  N_{0}^1 l(0)_1^2\rho^2 \vec{q}_{0}^1, \label{eqn:FDEC_equation}
\end{equation}
which needs to be mapped to the first order ellipsoid equation:
\begin{equation}
  \tilde{\vec{f}}(\rho, \phi) = R \cdot
  \begin{pmatrix}
    a \hat{x} \\
    b \hat{y} \\
    c\rho^2 /2 
  \end{pmatrix}.                              
\end{equation}
The rotation matrix $R$ reflects that the coordinate system may not be aligned with the axes of the fitted ellipsoid.  
The traditional principal component analysis (PCA) solves this by searching the maximal lengths of vectors from the
plane tangent to $\vec{\hat{f}}$ which must have the following form:

\begin{equation*}
  \hat{\vec{g}}(\rho, \phi) = A^T \cdot \left(
    \begin{matrix}
      \hat{x} \\
      \hat{y} 
    \end{matrix}
  \right).
  \end{equation*}
The length of any such vector is then:
\begin{equation*}
  \hat{\vec{g}}(\rho, \phi)^T \cdot \hat{\vec{g}}(\rho, \phi) =
\left(
    \begin{matrix}
      \hat{x} \\
      \hat{y}
    \end{matrix}
  \right) ^T  A \cdot 
  A^T \cdot \left(
    \begin{matrix}
      \hat{x} \\
      \hat{y} 
    \end{matrix}
  \right).
  \end{equation*}
  For a unit $\rho = \sqrt{\hat{x}^2 + \hat{y}^2} = 1$, such a distance will be maximal (resp. minimal) on the major (resp. minor) axes of the ellipsoidal belonging to the tangent plane of the surface. Therefore, $a^2$ and $b^2$ are the eigenvalues of the matrix $A^T A$. The associated eigenvectors $\hat{v}_a, \hat{v}_b$ allow to extract the principal directions of the ellipsoid belonging to the tangent plane: $v_a = A^T \hat{v}_a$ and $v_b = A^T \hat{v}_b$. This approach is similar to the solution method proposed in previous work \cite{REF:1, shaqfa2021b}, only with a different eigenproblem.

  If the expansion coefficients are expressed in a coordinate system aligned with $[v_a, v_b, v_a \times v_b]$ the following expression will hold:

  \begin{equation}
    \hat{\vec f}(\rho, \phi) =
    \begin{pmatrix}
      a & 0 \\
      0 & b
    \end{pmatrix}
    \cdot \left(
      \begin{matrix}
        \hat{x} \\
        \hat{y} 
      \end{matrix}
    \right) +  N_{0}^1 l(0)_1^2\rho^2 \vec{q}_{0}^1
    \equiv   \tilde{\vec{f}}(\rho, \phi) = 
    \begin{pmatrix}
      a & 0 \\
      0 & b
    \end{pmatrix}
    \cdot \left(
      \begin{matrix}
        \hat{x} \\
        \hat{y} 
      \end{matrix}
      \right)
      + \frac{c\rho^2}{2} \vec{e}_z,
  \end{equation}
  where the rotation matrix $R$ became identity. Assuming that $\vec{q}_0^1$ is normal to the tangent plane (i.e. parallel to $e_z$ in this frame), we can
  finally extract the parameter:
  \begin{equation}
    c = 2N_0^1l(0)_1^2 \mid\mid \vec{q}_0^1 \mid \mid.
  \end{equation}

\par
We generated caps with various $\theta_c$s from a unit sphere to test the accuracy of the computed curvature we obtain from FDEC by either the eigenvalue estimation or the size of the oriented bounding-box (OBB) method. For this, we reconstruct the surface using either k=1 or k=5. k=1 corresponds to the FDEC size while k=5 considers some higher modes. We then fit a bounding box around the reconstructed surface to obtain a, b and c. Such caps, of unit spheres, have a constant Gaussian curvature of $1$. From knowing $a$ and $c$ we can compute the constant Gaussian curvature of a spherical cap by the following formula:
\begin{equation}
    \kappa = \frac{4c^2}{(a^2+c^2)^2}.
\end{equation}
Table \ref{TAB:FDEC_size_comp} shows the results of size and curvature estimated by different methods. To benchmark the estimated cap size we used $a_{avg}$ instead of $a$ and it is computed as the average lengths of the major and intermediate eigenvalues ($a$ and $b$) of the reconstructed cap data. The first column contains the size of the input cap; the second, third and fourth column are the cap sizes obtained with the three methods described above (Eigenvalue, OBB k=1, OBB k=5). The results show that the Eigenproblem and OBB k=1 underestimate the cap size while OBB with k=5 leads to very good estimates not only of the cap size. The error on the curvature is larger because the curvature is proportional $1/L^2$.  This shows that as  the reconstruction includes more harmonics, one can get a better estimation of the cap size and the curvature.
\begin{table}[!ht]
\centering
\begin{tabular}{l|lll|}
 \cline{2-4}
 & \multicolumn{3}{c|}{FDEC size $(a_{avg}, c, \kappa)$}                                                                                                                              \\ \hline
\multicolumn{1}{|l|}{$\theta_c$ $(a, c, \kappa)$} & \multicolumn{1}{l|}{Eigenproblem} & \multicolumn{1}{l|}{OBB at $k =1$} & \multicolumn{1}{l|}{OBB at $k =5$} 
        \\ \hline
\multicolumn{1}{|l|}{$5^\circ$ $(0.174, 0.0038, 1)$}                    & \multicolumn{1}{l|}{$(0.116, 0.0031, 3.378)$}        & \multicolumn{1}{l|}{$(0.149, 0.004, 2.065)$}              & \multicolumn{1}{l|}{$(0.171, 0.0039, 1.130)$}             \\ \hline
\multicolumn{1}{|l|}{$10^\circ$ $(0.347, 0.0152, 1)$}                   & \multicolumn{1}{l|}{$(0.232, 0.012, 3.114)$}        & \multicolumn{1}{l|}{$(0.298, 0.016, 2.031)$}              & \multicolumn{1}{l|}{$(0.340, 0.016, 1.205)$}             \\ \hline
\multicolumn{1}{|l|}{$20^\circ$ $(0.685, 0.0602, 1)$}                   & \multicolumn{1}{l|}{$(0.505, 0.053, 2.536)$}        & \multicolumn{1}{l|}{$(0.588, 0.0635, 1.971)$}              & \multicolumn{1}{l|}{$(0.670, 0.0629, 1.172)$}          \\ \hline
\multicolumn{1}{|l|}{$50^\circ$ $(1.530, 0.357, 1)$}                   & \multicolumn{1}{l|}{$(1.060, 0.291, 2.534)$}        & \multicolumn{1}{l|}{$(1.350, 0.376, 1.587)$}              & \multicolumn{1}{l|}{$(1.500, 0.372, 1.127)$}           \\ \hline
\end{tabular}
\caption{Comparing FDEC sizes and curvatures from the estimated eigenproblem and OBB.}
\label{TAB:FDEC_size_comp}
\end{table}


%% file: app_A_fract_PSD_analytical.tex
\par
Out of the many types of fractal surfaces, engineered and real surfaces are mainly found to be self-similar or, more generally, self-affine \cite{Goedecke2013}. A self-similar surface extended in the Cartesian coordinates $f(x,y)$ is a surface that when magnified with a factor $\lambda$ is found to repeat itself such that $f(x,y) = \lambda f(x/\lambda, y/\lambda)$. Self-affine surfaces are surfaces that when magnified, scale as $f(x,y) = \lambda^{H} f(x/\lambda, y/\lambda)$, with $H$ representing the Hurst exponent (see \cite{Persson2005} for details). The corresponding fractal dimension is related to the Hurst exponent by $D=3-H$ (see \cite{Russ1994}).
\par
This paper is limited to isotropic self-affine surfaces where the fractal dimension is radially symmetric along with any open surface patch. This section will first briefly describe the derivation of Bessel functions of the first kind and zeroth order starting from the double Fourier transform; for details, we refer the readers to more detailed explanations elsewhere \cite{Fung1967, hecht_2002}. Then, we proceed to find the analytical correlation between the power spectrum computed using DH of a self-affine surface with the Hurst exponent. In other words, we are looking for the power decay in the power spectrum when a surface is decomposed with Fourier-Bessel basis functions.
\par
We know that the unnormalised forward double-Fourier transform of a continuous signal in the time domain can be written as:
\begin{equation}
    F(k_x, k_y) = \iint \limits_{-\infty}^{{~~\infty}} f(x,y) e^{-i (k_x x + k_y y)} \,d x \,d y,
\end{equation}
and the inverse transform can be written as:
\begin{equation}
    f(x, y) = \iint \limits_{-\infty}^{{~~\infty}} F(k_x, k_y) e^{i (k_x x + k_y y)} \,d k_x \,d k_y,
\end{equation}
where $k_x$ and $k_y$ are the integer frequencies (harmonics) of the expansion along the $x$ and $y$ axes.
\par
Assume $f(x,y)$ is a height map of a surface that is distributed onto a unit disk and parameterised as $f(\rho,\phi)$ with $k_{m} = \sqrt{k_x^2 + k_y^2}$ is the radial frequency, and thus it follows that $k_x = k_{m} \cos m$, $k_y = k_{m} \sin m$, $x = \rho \cos \phi$, $y = \rho \sin \phi$ with a finite area $\,d x \,d y = \rho \,d \rho \,d \phi$. Then by using the change of variables the latter equation becomes:
\begin{equation}
    f(\rho, \phi) = \iint \limits_{-\infty}^{{~~\infty}} F(k_m, m) e^{i(k_m \rho \cos \phi \cos m + k_m \rho \sin \phi \sin m)} \rho \,d \rho \,d \phi.
\end{equation}
Because we know that $\cos (\phi - m) = \cos \phi \cos m + \sin \phi \sin m$, we have:
\begin{equation}
    f(x, y) = \iint \limits_{-\infty}^{{~~\infty}} F(k_m, m) e^{i k_m \rho \cos(\phi-m)} \rho \,d \rho \,d \phi.
\end{equation}
Then the forward transform over a unit disk can also be written as:
\begin{equation}
    F(k_m, m) = \int_{0}^{1} \int_{0}^{2\pi} \rho f(\rho, \phi) e^{-i k_m \rho \cos(\phi-m)}  \,d \phi \,d \rho.
    \label{EQN:Bessel_FOUR_PSD}
\end{equation}
And from what we know, the Bessel function of the first kind for the zeroth order---an even function---can be alternatively written as \cite{REF:2, Watson1944} (notice the exponent sign change):
\begin{equation}
    J_0 (u) = J_0 (-u) = \frac{1}{2\pi} \int_{0}^{2\pi} e^{i u \cos(v)}  \,d v.
    \label{EQN:Bessel_ALternative_int}
\end{equation}
\par
The exclusive use of the first kind Bessel function with $m = 0$ can be limited to describe only functions that are radially symmetric (see Fig.~\ref{FIG:BESSEL_FOURIER_BASIS} where $m=0$ for an illustration). Now, we assume that the map $f(\rho, \phi)$---though not the actual height map of the surface---is \say{statistically} isotropic onto a unit disk. Then, $f(\rho, \phi)$ is radially symmetric about the origin of the disk so that it can be infused into $f(\rho)$. This directly implies that the Fourier-Bessel functions of order $m = 0$ will be sufficient for decomposing isotropic self-affine surfaces. By substituting Eq.~(\ref{EQN:Bessel_ALternative_int}) into Eq.~(\ref{EQN:Bessel_FOUR_PSD}), we obtain:
\begin{equation}
    F(k_0) = 2\pi \int_{0}^{1} \rho f(\rho) J_0(k_{0} \rho) \,d \rho.
    \label{EQN:PSD_int}
\end{equation}
\par
Now, to define the linkage between the PSD computed from the DH basis functions and the Hurst exponent $H$, we depend on the extension of the Wiener-Khintchine theorem for autocorrelation functions by Fung (1967) \cite{Fung1967}. In that paper, the author extended the one-dimensional Fourier analysis to a multidimensional analysis and discovered that it can be extended by Bessel functions of the first kind and with order $m = 0$, similar to the assumptions we have reached.
\par
From the Wiener-Khintchine theorem, we can write the height-height autocorrelation function (ACF) as a function of the power spectrum. Let $\Psi(\rho) = \langle f(\rho) f(0) \rangle$ be a height-height ACF, where $\langle\dots\rangle$ is the ensemble averaging of the height maps with a spatial delay; this is not to be confused with the inner dot product $\langle \cdot, ~\cdot \rangle$ that was assigned previously. Ensemble averaging is used to average the heights over a set of independently generated surfaces that have a common statistical property. The power spectrum can then be written as:
\begin{equation}
    ||F(k_0)||^{2} = \lim_{R\to\infty} \frac{2}{R^2} \int_{0}^{R}  \langle f(\rho) f(0) \rangle \rho J_0(k_0 \rho) \,d \rho,
    \label{EQN:PSD_W-K}
\end{equation}
where $R$ is the radius of the considered surface in this analysis.
\par
From Eq.~(\ref{EQN:PSD_W-K}), we are interested in the limiting behaviour of Bessel functions on an infinitely extended disk where the PSD implies $k_0 R$ in Eq.~(\ref{EQN:PSD_W-K}) or the equivalent to $u \to \infty$ as in Eq.~(\ref{EQN:Bessel_ALternative_int}). Then, the Bessel function of the first kind can be replaced by the asymptotic form in Eq.~(\ref{EQN:Bessel_ASYM}) and plugged into Eq.~(\ref{EQN:PSD_int}) to obtain:
\begin{equation}
    ||F(k_0)||^2 \approx \frac{1}{R^2}\sqrt{\frac{8}{\pi k_0}} \int_{0}^{R} \langle f(\rho) f(0) \rangle \sqrt{\rho} \sin\Bigg(k_0 \rho+\frac {\pi }{4}\Bigg) \,d \rho.
    \label{EQN:PSD_int2}
\end{equation}
With a properly decaying power-law, the latter formula can be used for generating rough fractal surfaces if we randomise the phase shift in the sine term using a uniformly distributed noise.
\par
Now, we can re-write Eq.~(\ref{EQN:PSD_int2}) as if we zoom in to the disk details with a factor of $\lambda$, then $\rho = \rho^\prime /\lambda \to \,d \rho = \,d \rho^\prime/\lambda$ and let $k_0 \rho + \pi/4 \approx \lambda \rho$ (swapped for simplicity as the asymptotic roots approximately equal the actual ones when $\lambda \to \infty$) as well as assuming a unit disk ($R=1$) to ignore the size normalisation:
\begin{equation}
    ||F(\lambda)||^2 \approx \sqrt{\frac{8}{\pi k_0}} \int_{0}^{1} \lambda^{-3/2} \langle f(\rho^\prime/\lambda) f(0) \rangle \sqrt{\rho^\prime} \sin\big(\rho^\prime\big) \,d \rho^\prime.
\end{equation}
Then, from the definition of a self-affine surface, we can define $\langle \lambda^{2H} f(\rho^\prime/\lambda) f(0) \rangle = \langle f(\rho^\prime) f(0) \rangle$. The power spectrum then becomes:
\begin{equation}
    ||F(\lambda)||^2 \approx \sqrt{\frac{8}{\pi k_0}} \int_{0}^{1} \lambda^{-3/2} \lambda^{-2H} \langle f(\rho^\prime) f(0) \rangle \sqrt{\rho^\prime} \sin\big(\rho^\prime\big) \,d \rho^\prime.
    \label{EQN:PSD_H_Bessel_pre}
\end{equation}
In the herein proposed analysis approach, we only evaluate this ACF function $\Psi(\rho)$ at the eigenvalues of the Fourier-Bessel functions and not the eigenvalues of the Fourier transform ($k_0$). Let us assume that $\lambda$ is the eigenvalues (zeros) of the Bessel functions of the first kind and zeroth order such as $\lambda = l(0)_k$. Thus, wherein $\lambda\rho^{\prime} \to \infty$, Eq.~(\ref{EQN:PSD_H_Bessel_pre}) can be written as:
\begin{equation}
    ||F(\lambda)||^2 \approx \sqrt{\frac{8}{\pi k_0}}~\Big(\lambda^{-2 (3/4+H)}\Big) \int_{0}^{1} \langle f(\rho^\prime) f(0) \rangle \sqrt{\rho^\prime} \sin\big(\rho^\prime\big) \,d \rho^\prime.
    \label{EQN:PSD_H_Bessel}
\end{equation}
From the latter, we can see that an isotropic self-affine surface will scale as:
\begin{equation} 
    ||F(\lambda)||^2 \propto \lambda^{-2 (3/4+H)},
    \label{EQN:summary_PSD_vs_H}
\end{equation}
or in terms of the fractal dimension ($D$), as:
\begin{equation}
    ||F(\lambda)||^2 \propto \lambda^{-(15/2-2D)}.
    \label{EQN:summary_PSD_vs_D}
\end{equation}
\par
Equation~(\ref{EQN:PSD_H_Bessel}), shows the power decay of the power spectrum that occurs with increasing eigenvalues of the Bessel functions. Note that in this paper, the power spectrum means the moments computed from the zeroth-order basis only as $\text{PSD}_{\text{DH}, m = 0}$ [as the ACF implies in Eq.~(\ref{EQN:PSD_W-K})]. Experimentally, this should be realised from plotting the computed PSD using Fourier-Bessel basis functions versus the eigenvalues (wave number) on a log-log scale graph by fitting a power law function (with the least squares algorithm).

\begin{figure*}[!ht]
    \centering
    \includegraphics[width=0.95\textwidth]{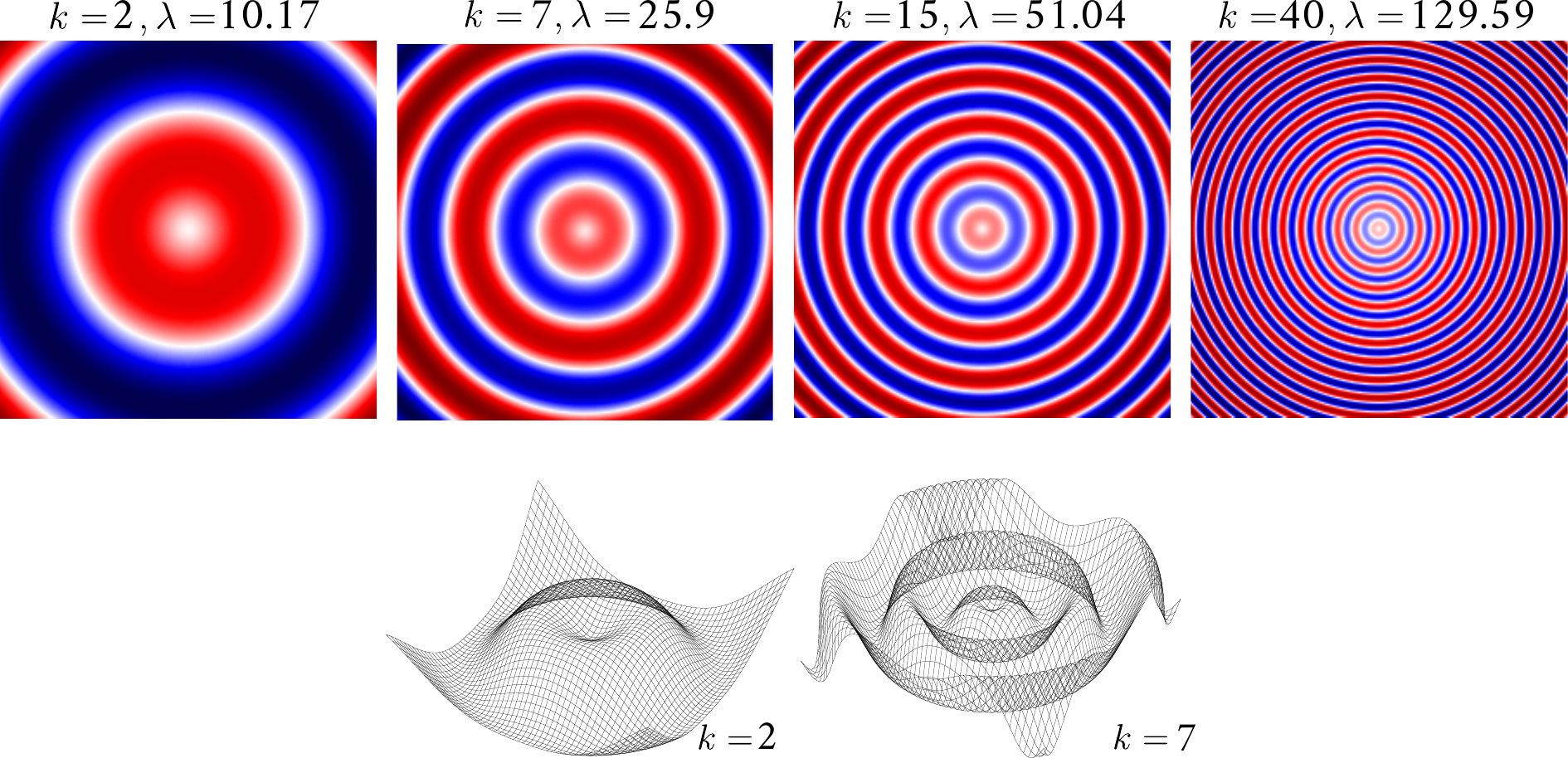}
    \caption{An approximate form of the Bessel functions of order $m=0$ on a rectangular grid $f(\rho, \phi) = \sqrt{\rho}\sin(\lambda \rho)$ (cf. with Figure \ref{FIG:BESSEL_FOURIER_BASIS_3D} for $m=0$).
    \label{FIG:SIMPLIFIED_BESSEL}}
\end{figure*}

A brief commentary on the results of the derivation of Eq.~(\ref{EQN:summary_PSD_vs_H}):
\begin{itemize}
    \item Comment 1: In this analysis, we are not concerned about the topological features and the actual height map. Instead, we focus only on how these maps scale at different eigenvalues, which is the statistical property comprising such surfaces. Thus, we consider $f(\rho)$ rather than $f(\rho, \phi)$. This does not mean that our height map is actually radially symmetric, but rather the surface is isotopically self-affine (in the power law).
    
    \item Comment 2: Integrating the harmonics of $m=0$ only comes as a natural extension of the Wiener-Khintchine theorem \cite{Fung1967} for expanding ACF with 1D Fourier filters into a multidimensional analysis. Also, the zeroth-order Bessel function occurs for all analysis degrees from $k=0$ up to any arbitrary $k$, but for other orders, such as the PSD when $m = 5$ for example, this is repeated for degrees $k \geq 5$ only and cannot be found for lower degrees or larger wavelengths. This means that not all orders are defined for all scales, apart from $m=0$.
    
    \item Comment 3: Analysing randomly generated surfaces by only randomising the phase using Fourier filters will result in a non-smooth (noisy) decay because the interharmonics are fit based on the Bessel eigenvalues rather than the actual Fourier harmonics (eigenvalues) that were used to generate the surface. Some interharmonics can also be poorly correlated, leading them to practically hit a zero, which was defined as below $10^{-8}$ in our experiments; thus we can exclude these values from the analysis.
    
    
    
    \item Comment 4: Analysing surfaces with the ACF in Eqs.~(\ref{EQN:PSD_int2})--(\ref{EQN:PSD_H_Bessel}) while using $\lambda = k_0$ (assuming Fourier basis eigenvalues) will increase the spectrum decay by a factor of $\lambda^{-1/2}$ due to the constant multiplier $\sqrt{\frac{8}{\pi k_0}}$ in the same equation. This makes the overall power decay of $\lambda^{-2(1+H)}$, which corresponds to the result found in Persson et al. (2005) \cite{Persson2005}.
    
    \item Comment 5: The simplified function extracted from Eq.~(\ref{EQN:PSD_H_Bessel_pre}), $f(\rho, \phi) = \sqrt{\rho}\sin(\lambda \rho)$, simply resembles a wave function that travels from the origin of a disk, and the number of oscillations is related to the index of the eigenvalue of the Bessel functions $(k)$ (see Fig.~\ref{FIG:SIMPLIFIED_BESSEL}). With a simple random phase shift on the sine term, one can use Eq.~(\ref{EQN:PSD_int2}) for generating fractal surfaces with a proper power law.
    
\end{itemize}

%% file: app_A_generating_self_affines.tex
\par
Because generating self-affine scaling (fractal) surfaces is important for many theoretical and practical applications, several papers and methods have been proposed for handling this. For the completeness of the paper, we herein briefly summarise the creation of an isotropic power-law for the FFT filter method. For more details about the generation process, we recommend reading Persson et al. (2005) \cite{Persson2005} and Jacobs et al. (2017) \cite{Jacobs2017}. The herein used terminology, notation and scaling relations are derived by Persson et al. \cite{Persson2005}.
\par
Persson et al. used an approach similar to the one in \ref{APP:fract_PSD} to find that the PSD is proportional to the Hurst exponent $H$ of an isotropic surface as:
\begin{equation}
    C(q)^{iso} \propto q^{-2 (1+H)},
    \label{EQN:summary_PSD_vs_H_Persson}
\end{equation}
where $q$ is the wavevector ($q = k_m$ in this paper, though this section is kept consistent with the original notations). The largest wavevector (smallest wavelength) is typically called the cut-off wavevector $q_s$ at the atomic spacing, wherein this measurement is mostly limited by the measurement accuracy. However, the smallest wavevectors (largest wavelength) is $q_l$. Normally for $q<q_l$, the power law is considered constant and depends on the surface treatment; the constant regime is extended to the roll-off wavevector $q_r$. In general, $q_r<q_l<q_s$ for any assumed power law of an isotropic self-affine scaling surface, and the power-law is constituted by knowing $H$ as:
\begin{equation}
C(q)^{iso} = 
\left\{
\begin{array}{ll}
      0 & q < q_r\\
      q_l^{-2(1+H)}=const. & q_r \leq q \leq q_l\\
      q^{-2(1+H)} & q_l \leq q < q_s\\
      0 & q \geq q_s\\
\end{array}
\right..
\label{EQN:iso_power_law_PSD}
\end{equation}
An example of such a power law is shown in Fig.~\ref{FIG:PSD_SEFL_AFFINE} and used to generate self-affine surfaces for the validation section in this paper.

%% file: app_A_roughness_projection.tex
\par
We project rough fractal disks onto spherical caps to generate rough surfaces with a predefined curvature and a power law. Smooth spherical caps have a constant Gaussian curvature everywhere on the surface, and the principal curvatures can be computed as $\kappa = 1/R^2_{\theta_c}$, where $R_{\theta_c}$ is the radius of the sphere or cap. We can project the generated fractal disks onto a predefined spherical cap using a proper mapping method, for which we herein used the inverse Lambert azimuthal equal-area projection as in our previous work \cite{shaqfa2021b}.
\par
More specifically, the Lambert's projection $\tau_{l}(X,Y,Z)$ and the Lambert's inverse projection $\tau^{-1}_{l}(x, y)$ are defined as:
\begin{equation}
\label{EQN:LAMBERT_PROJ}
\begin{split}
    \tau_{l}(X,Y,Z) &= (x,y) =  \left(\sqrt{\frac{2}{1-Z}} X, \sqrt{\frac{2}{1-Z}} Y \right), \\
        \tau^{-1}_{l}(x, y) &= (X,Y,Z) =
        \left(\sqrt{1-\frac{x^2 + y^2}{4}}x, \sqrt{1-\frac{x^2 + y^2}{4}}y,-1+\frac{x^2 + y^2}{2}\right).
\end{split}
\end{equation}
The Lambert's inverse projection $\tau^{-1}_{l}(x, y)$ can be applied over a rescaled disk denoted by $\mathcal{D}_l$ and its generated height-map $z = h(x,y)$ to get the needed spherical cap $\mathbb{S}^2_{~\theta \leq \theta_c}$. To project a unit disk onto a unit sphere with a prescribed half-angle $\theta_c$, we need to scale the radius of the disk with a scaling factor, derived in our previous work \cite{shaqfa2021b}, to achieve the required half-angle $\theta_c$:
\begin{equation}
    \label{EQN:SCALE_LAMBERT}
    r_l = \sqrt{2(1-\cos\theta_c)}.
\end{equation}
To this end, we can project a smooth circular disk onto a smooth spherical cap to define the corresponding domain. Next, we also need to project the height map (generated surface roughness over the disk domain) onto the smooth spherical cap $\mathbb{S}^2_{~\theta \leq \theta_c}$ such that we preserve the normals of the map (normal-preserving map). Preserving normals means that the normals of asperities, relative to their vicinity, must be pointing in the same direction as they did on the original disk. The simplest projection that can preserve normals on caps is the radial projection $\tilde{\tau}_{h}$, which rescales every point by $(1 + \sqrt{s_c} h(x,y))$. The final coordinates of the height maps are $(\tilde{X},\tilde{Y},\tilde{Z})$:
\begin{equation}
\tilde{\tau}_{h} \circ \tau^{-1}_{l} = 
\left\{
\begin{array}{l}
      \tilde{X} = X \big(1 + \sqrt{s_c} h(x,y)\big)\\
      \tilde{Y} = Y \big(1 + \sqrt{s_c} h(x,y)\big)\\
      \tilde{Z} = Z \big(1 + \sqrt{s_c} h(x,y)\big)\\
\end{array}
\right.,
\label{EQN:final_proj_coord}
\end{equation}
where $s_c$ is the ratio between the flat configuration of a spherical cap (prescribed with $\theta_c$ and radius $R_{\theta_c}$) to the smooth (nominal) area of the generated rough disk.
\par
\begin{figure}
    \centering
    \includegraphics[width=1.0\linewidth]{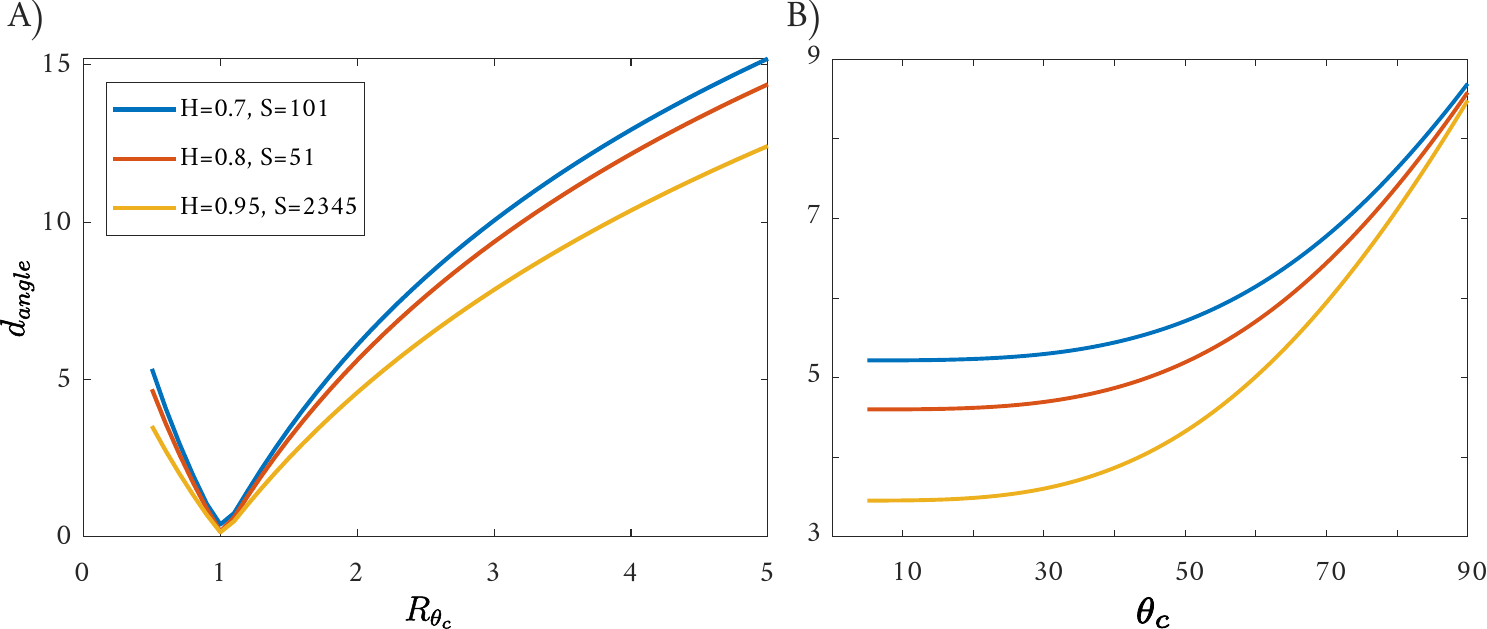}
    \caption{Angular distortion associated with the roughness projection. A) Angular distortion computed in degrees, for some of the selected patches with $\theta_c = 10^\circ$ and various seed numbers (S) when changing the radii $R_{\theta_c}$. B) Angle distortion when changing the half-angle $\theta_c$.}
    \label{FIG:angle_distortion}
\end{figure}

This projection is not an isometric projection as it does not preserve neither the angles nor the areas of the elements. The projection is thus a source of angular distortion that could distort the PSD as it varies from the original smooth configuration to the newly introduced cap configuration. Also, the final radially projected height map is not, in general, area-preserving. We propose to compute the angle distortion as an average of the absolute difference between the disk and cap meshes per element to be measured, as in our previous work \cite{shaqfa2021b}:
\begin{equation}
    d_{angle} = \underset{[v_i,v_j,v_k]}{\text{mean}} {|\angle([f(v_i),f(v_j),f(v_k)]) - \angle([v_i,v_j,v_k])|}.
    \label{EQN:angle_distortion}
\end{equation}
Figure \ref{FIG:angle_distortion}A shows $d_{angle}$ for three surface patches with $\theta_c = 10^\circ$ using the inverse Lambert azimuthal projection. The optimal point for $d_{angle}$ is when $R_{\theta_c} = 0.35$. Figure \ref{FIG:angle_distortion}B shows the $d_{angle}$ dependency on of the cap's half-angle $\theta_c$ and it shows that with small $\theta_c < 20^\circ$ the distortion is minimal. In general, when $R_{\theta_c} > 1$ and ${\theta_c} > 20^\circ$, there is a large angular distortion in the curved rough surface. Also, the higher the $H$ of the surface, the lower the obtained distortion as the surface texture becomes smoother.

%% file: app_B_boundaryconditions.tex
\par
\begin{figure*}[!ht]
    \centering
    \includegraphics[width=1.0\textwidth]{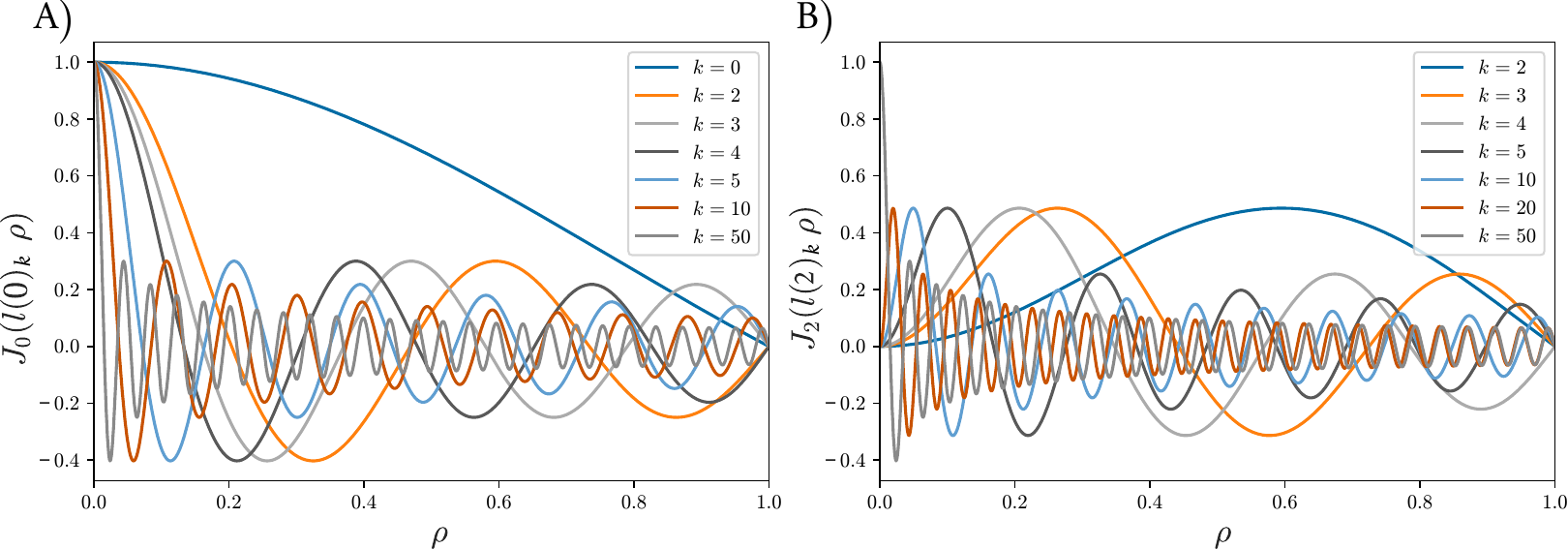}
    \caption{The unnormalized Bessel functions with a Dirichlet BC for $\rho \in [0, 1]$ as the functions change signs $k-m$ times over the orthogonality domain.
    \label{FIG:BESSEL_Dirichlet_BC}}
\end{figure*}

Similar to Fig.~\ref{FIG:BESSEL} with the same orders and indices, Fig.~\ref{FIG:BESSEL_Dirichlet_BC} shows Bessel functions, though with a Dirichlet BC imposed at $\rho(\phi) = 1, \forall \phi \in [0,2\pi]$ with $m\in\{0,2\}$ and different $k$ indices. As can be seen, the end points all collapse at zero, which will obviously not retrieve the right reconstruction of the free edge on the disk.
\par
\begin{figure*}[!ht]
    \centering
    \includegraphics[width=1.0\textwidth]{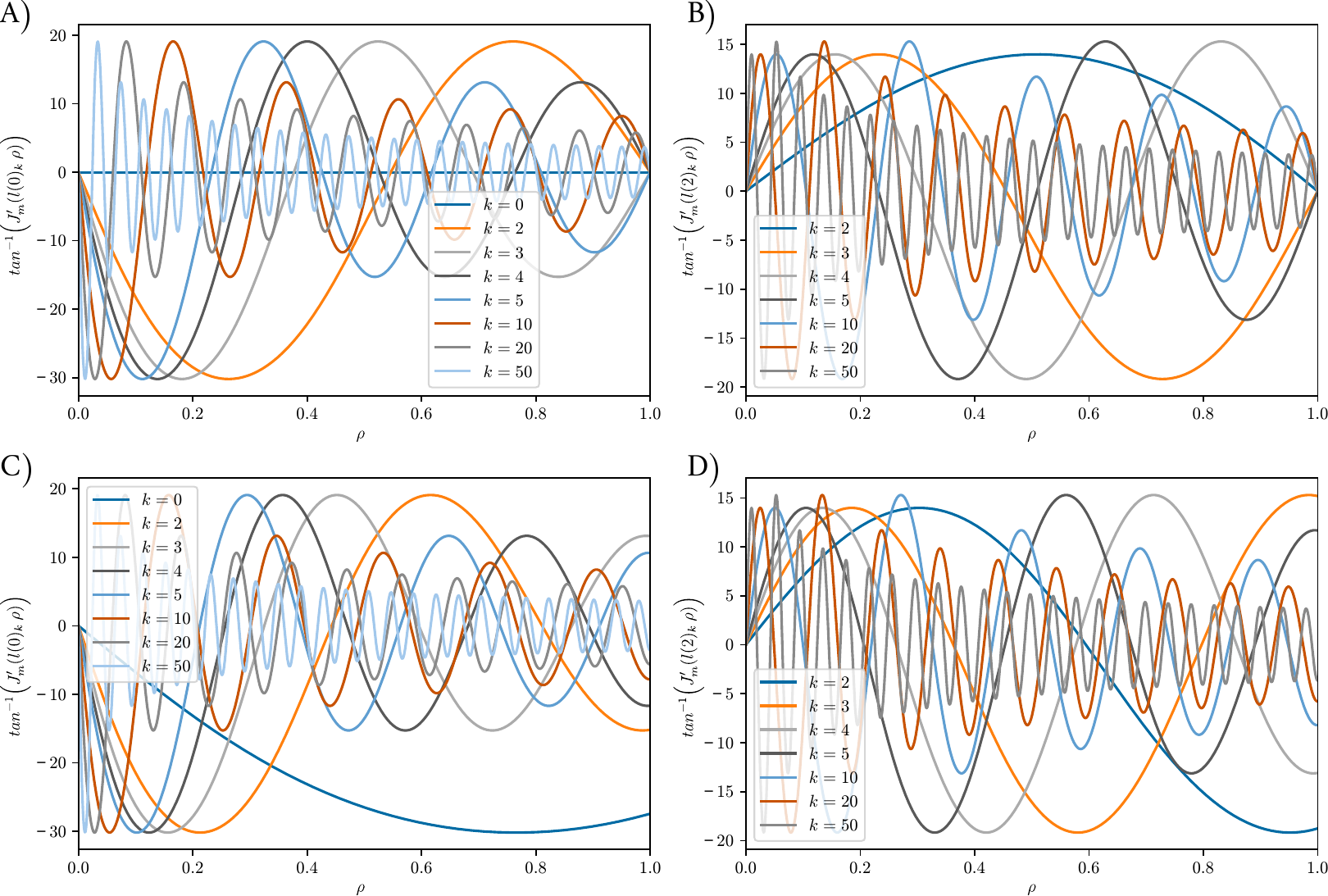}
    \caption{The slopes of Bessel basis functions in degrees. (A), (B) The slopes of the plotted basis functions in Figure \ref{FIG:BESSEL} corresponded directly, one-to-one, with the Neumann BC as ${\partial J_m (l(m)_k \cdot 1)}/{\partial \rho} = 0$. (C), (D) The slopes shown in Figure \ref{FIG:BESSEL_Dirichlet_BC} corresponded to the Dirichlet BC as $J_m (l(m)_k \cdot 1)
    = 0$.
    \label{FIG:BESSEL_BC_SLOPE}}
\end{figure*}

Instead, we applied the Neumann BC, and here we show how the imposed BC influences the basis functions. Figure \ref{FIG:BESSEL_BC_SLOPE} shows the change in slopes (in degrees) of both Bessel function with the Neumann BC as shown in insets (A) and (B), and with the Dirichlet BC in (C) and (D). The zero slope for the Neumann BC only affects the endpoint $\rho = 1$, and the slopes around that vicinity change from one basis function to another, which indeed will arbitrarily change with the realised expansion coefficients and will depend on the mesh size of the used finite grid for the analysis and reconstruction stages.

%% file: app_B_comparison_DHA_vs_SCHA.tex
\par
\begin{figure*}
\centering
\includegraphics[width=0.9\textwidth]{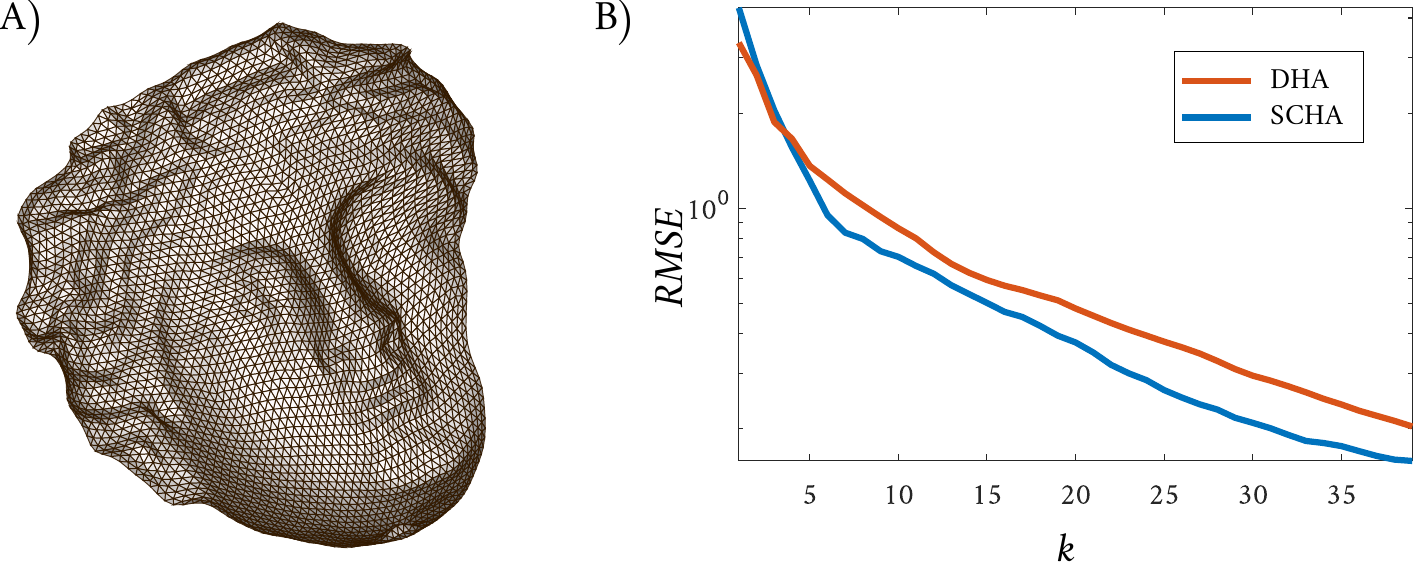}%
\caption{A comparison of the analysis and reconstruction between the SCHA method in Shaqfa et al. (2021b) \cite{shaqfa2021b} and the DHA approach proposed by this paper. (A) The DHA-reconstructed surface mesh at $k = 40$. (B) The convergence of the root mean square error (RMSE) computed for both SCHA and DHA against the original input surface.
\label{FIG:DHA_vs_SCHA}
}
\end{figure*}

We analysed and reconstructed the 3D benchmark face used in our SCHA paper \cite{shaqfa2021b}, with Fig.~\ref{FIG:DHA_vs_SCHA} showing that the SCHA approach converges better than the DHA. However, such results are highly dependent on the parameterisation method that is used, and as these both differ, we cannot directly conclude which method is superior in spatial convergence. The normalised RMSE computed for DHA was $0.203526$, while from SCHA we obtained $0.158282$ as $k=40$, which is a difference that cannot be spotted by an unaided eye; notice the semi-log scale used in Fig.~\ref{FIG:DHA_vs_SCHA}(B).
\par
The added value of this approach is indeed more apparent when we compare the analysis and reconstruction times of both approaches. Here, the DHA took well under half a minute, while the SCHA took more than $16$ hours to construct the basis functions only due to the poor convergence of the hypergeometric functions; see \cite{shaqfa2021b} for details. In this context, we consider the DHA approach as far more advantageous than the SCHA, and the DHA method can reach higher degrees that correspond to a better resolution.
\par
\begin{figure*}
\centering
\includegraphics[width=0.80\textwidth]{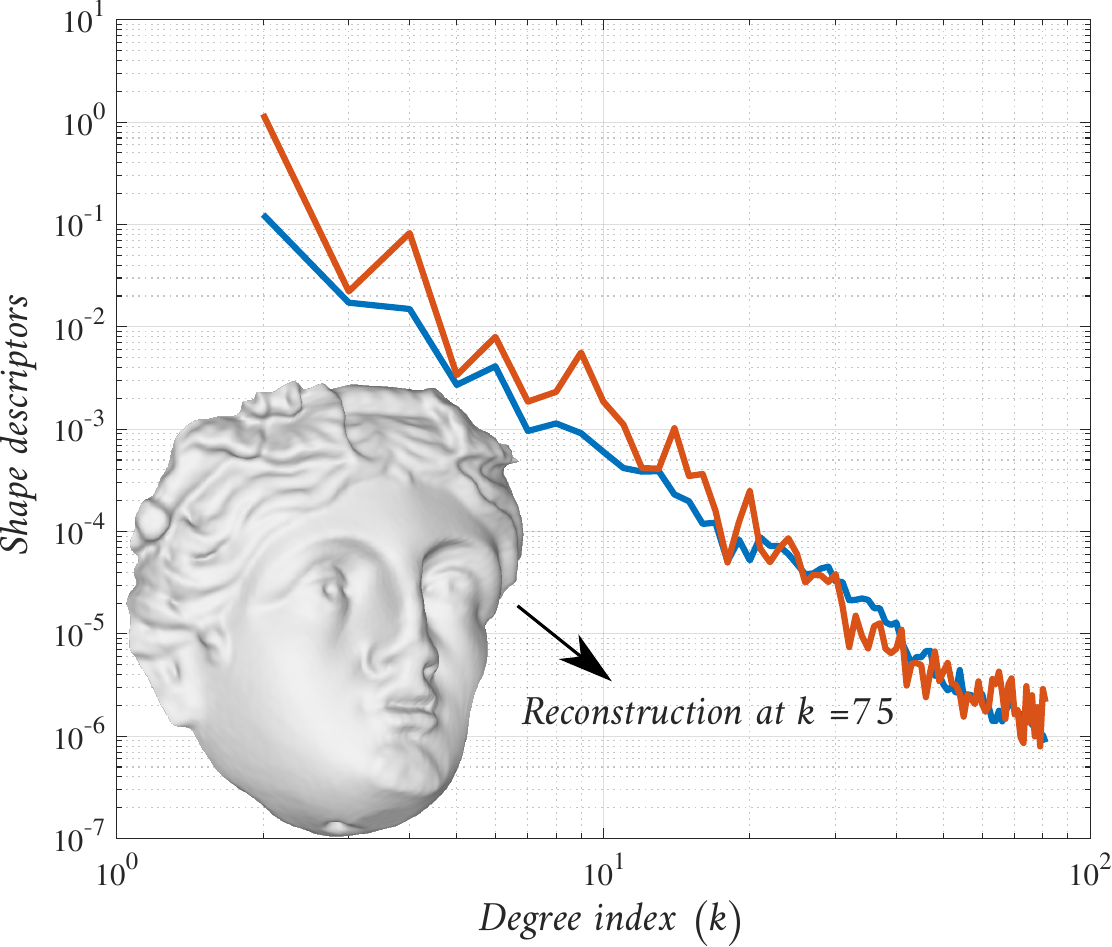}%
\caption{The reconstruction of the 3D face at $k = 75$. The solid red line shows only the normalised descriptors of the zeroth-order harmonics, while the solid blue line shows the normalised descriptors for all harmonics, including all orders up to $k = 75$. The computed RMSE at that degree was $0.077543$.
\label{FIG:rec_3D_face_k_75}
}
\end{figure*}

Using DHA, we constructed the results of the face with $k=75$ and a finer mesh size. In this case, the resulting normalised RMSE was only $0.077543$. The time needed for this analysis was less than $5$ minutes, which was accomplished by exploiting a parallel least-squares solver on an Intel Core i7-11370H CPU with 4 cores and 8 threads and a maximum single-core frequency of 4.80 GHz. Figure \ref{FIG:rec_3D_face_k_75} depicts the reconstruction results as well as the computed normalised shape descriptors for all orders as well as for the zeroth order. It should be mentioned that from Fig.~\ref{FIG:rec_3D_face_k_75}, the problem was analysed for $k=80$, though after $k=75$, the model was over-fitted as the PSD began to oscillate, and more data points were required to fit the shorter wavelengths.
